\documentclass[a4paper,openright,12pt]{article}
\usepackage{geometry}
\usepackage{amsmath}
\usepackage{amsfonts}
\usepackage{amssymb}
\usepackage{graphicx}
\usepackage{verbatim}
\immediate\write18{texcount -tex -sum  \jobname.tex > \jobname.wordcount.tex}
\usepackage[numbers]{natbib}
\newcommand{\Sref}{Section \ref}
\newcommand{\eref}{\eqref}
\newcommand{\Fref}{Figure \ref}
\newcommand{\Tref}{Table \ref}
\def\tr{{\rm{tr}}}
\usepackage{hyperref}
\providecommand{\keywords}[1]
{
  \small	
  \textbf{\textit{Keywords---}} #1
}
\title{Turing patterns in a diffusive Holling--Tanner predator-prey model with an alternative food source for the predator}
\author{Claudio Arancibia--Ibarra$^{1,2}$, Michael Bode$^{1}$, Jos\'e Flores$^{3}$,\\
Graeme Pettet$^{1}$  and Peter van Heijster$^{1}$\\
\small $^{1}$School of Mathematical Sciences, Queensland University of Technology (QUT), \\
\small Brisbane, Australia. \\
\small $^{2}$Facultad de Educaci\'on, Universidad de Las Am\'ericas (UDLA), \\
\small Santiago, Chile. \\
\small $^{3}$Department of Mathematics, The University of South Dakota (USD), \\
\small Vermillion, South Dakota, USA \\
\small Email: claudio.arancibia@hdr.qut.edu.au,  michael.bode@qut.edu.au, Jose.Flores@usd.edu, \\
\small graeme.pettet@alumni.qut.edu.au, petrus.vanheijster@qut.edu.au. 
}

\begin{document}
\maketitle

\begin{abstract}
In this manuscript, we consider temporal and spatio-temporal modified Holling--Tanner predator-prey models with predator-prey growth rate as a logistic type, Holling type II functional response and alternative food sources for the predator. From our result of the temporal model, we identify regions in parameter space in which Turing instability in the spatio-temporal model is expected and we show numerical evidence where the Turing instability leads to spatio-temporal periodic solutions. Subsequently, we analyse these instabilities. We use simulations to illustrate the behaviour of both the temporal and spatio-temporal model.
\end{abstract} \hspace{10pt}


\section{Introduction}
One of the main problems in the ecological sciences is to understand the complex dynamical behaviour of the interaction between species. These interactions are becoming increasingly important in both ecology \cite{mondal, santos, turchin} and applied mathematics \cite{saez,yu,zhao}. The goals of the analysis of these interactions are to describe different behaviours between species, to understand their long term behaviour, and to predict how they respond to management interventions \cite{hooper,may2}. Dynamic complexities in such models, and in particular the Holling--Tanner predator-prey models, are of particular mathematical interest on both temporal \cite{arancibia7,arancibia3,arrows} and spatio-temporal domains \cite{banerjee1,banerjee2,ghazaryan}. 

The Holling--Tanner model has been used extensively to model many real-world predator-prey interactions \cite{turchin,hanski2,hanski,wollkind,hanski3,andersson,erlinge,hansson}. For instance, Hanski {\em et al.} \cite{hanski3} used the original Holling--Tanner model to investigate the multi-annual oscillation of field vole (\textit{Microtus agrestis}) in Fennoscandia. This oscillation is generated by the predator-prey interaction between the rodent and the least weasel (\textit{Mustela nivalis}) and the authors postulated that the least weasel population causes a delayed density dependence and therefore an oscillation phenomenon. However, the least weasel can switch its main food source depending on the proportion of the prey density. In particular, the weasel has three main food sources available birds and birds' eggs (\%5 of weasels diet), rabbit (\%25 of weasels diet) and small rodents (\%68 of weasels diet) \cite{mcdonald}. This characteristic was not considered in investigations associated with the oscillation of field vole in Fennoscandia which is affected by least weasel predators \cite{turchin,hanski2,hanski,wollkind,hanski3,andersson,erlinge,hansson}. That is, these studies do not consider that since the predator is a generalist they can survive under different environments and utilise a large range of food resources. Instead of adding more species to the model, we assume that these other food sources are abundantly available \cite{mcdonald} and model this characteristic by adding a positive constant $c$ to the environmental carrying capacity for the predator \cite{aziz}. Therefore, we have a modification to the prey-dependent logistic growth term in the predator equation, namely $K(N)=hN$ is replaced by $\overline{K}(N)=hN+c$. Additionally, there exists evidence that the amplitude of the oscillation in the predator-prey interaction is affected by geographic changes, since the predator density varies form north to south in Fennoscandia \cite{hanski}. Therefore, we also allow both species to diffuse. 

The model of interest is
\begin{equation}\label{eq01}
\begin{aligned}
\frac{\partial N}{\partial t} &= rN \left( 1-\frac{N}{K}\right)-\frac{qNP}{N+a}+D_1\triangledown^2 N, \\ 
\frac{\partial P}{\partial t} &= sP\left( 1\ -\frac{P}{hN+c}\right)+D_2\triangledown^2P.
\end{aligned}  
\end{equation}
In system \eref{eq01}, $N$ and $P$ indicate the prey and predator population sizes respectively, the predator and prey population contain logistic growth functions and the predator environmental carrying capacity is a prey dependant. Moreover, the functional response is hyperbolic in form and is referred to as a Holling Type II functional response \cite{turchin,may}. Additionally, $r$ and $s$ are the intrinsic growth rate for the prey and predator respectively, $h$ is a measure of the quality of the prey as food for the predator respectively, $K$ is the prey environmental carrying capacity, $q$ is the maximum predation rate per capita, $a$ is half of the saturated level and $c$ is considered the level of predators that are fed by the alternative food. We assume all parameters to be positive and, for ecological reasons, $a<K$. The predator and prey are assumed to diffuse through the spatial domain with diffusive coefficient $D_1$ for the prey population respectively $D_2$ for the predator population and $\nabla^2$ is the standard Laplacian operator.

The aim of this manuscript is to study the spatio-temporal dynamics of the modified Holling--Tanner predator-prey model \eref{eq01}. We will show that the addition of the alternative food source for the predator will lead to different Turing patterns when compared to the original Holling--Tanner model (i.e. $c=0$ in system \eref{eq01}) \cite{banerjee2}. While the original Holling--Tanner model is singular for $N=0$, \eref{eq01} is not singular and there exist system parameters that lead to spatio-temporal periodic solutions, see for instance \Fref{F05} in Subsection \ref{nusip2}. This is not possible for the original Holling--Tanner model \cite{banerjee2,ma}.

The temporal properties of the diffusion-free model were studied in \cite{arancibia,arancibia2,arancibia5,arancibia9} and are briefly discussed in \Sref{model}. In this section we also discuss the basins of attraction of the equilibrium points of the diffusion-free system. In \Sref{stem} we determine the Turing space where the Turing patterns occur and we present numerical simulations for different types of Turing patterns in one and two space dimension. Finally, in \Sref{con} we compare the Turing space to the model without alternative food as studied in \cite{banerjee2,ma} and we discuss the ecological implications.

\section{Temporal Model}\label{model}
In order to simplify the analysis we introduce dimensionless variables by setting $u:=N/K$, $v:=P/(hK)$, $S:=s/r$, $C:=c/(hK)$, $A:=a/K$, $Q:=qh/(Kr)$, $\tau:=rt$, $x=X\sqrt{r/D_1}$ and $d=D_2/D_1$. By substitution of these new parameters and variables into the one-dimensional version of \eref{eq01} we obtain
\begin{equation}\label{eq03}
\begin{aligned}
\frac{\partial u}{\partial\tau} &= uF(u,v)+\nabla^2u = u\left(\left(1-u\right)-\frac{Qv}{u+A}\right)+u_{xx}, \\
\frac{\partial v}{\partial\tau} &= vG(u,v)+d\nabla^2v = Sv\left(1-\frac{v}{u+C}\right)+dv_{xx}.
\end{aligned}
\end{equation}
System \eref{eq03} is defined in $u(x,t)\in\mathbb{R}_{\geq0}\times\mathbb{R}_{\geq0}$ and $v(x,t)\in\mathbb{R}_{\geq0}\times\mathbb{R}_{\geq0}$ and we first recall the stability of the equilibrium points of the diffusion-free system\footnote{Note that a different nondimensional version of the diffusion-free system \eref{eq01} was studied in \cite{arancibia,arancibia2,arancibia5,arancibia9}}
\begin{equation}\label{ODEeq03}
\begin{aligned}
\frac{du}{d\tau} = uF(u,v), \\
\frac{dv}{d\tau} = vG(u,v).
\end{aligned}
\end{equation}
System \eref{ODEeq03} is of Kolmogorov type, that is, the axes $u=0$ and $v=0$ are invariant and solutions curves initiated in the first quadrant (including the axes) stay in the first quadrant. The $u$ nullclines are $u=0$ and $v=(u+A)(1-u)/Q$, while the $v$ nullclines are $v=0$ and $v=u+C$. Hence, the equilibrium points for this system are  $(0,0)$, $(1,0)$, $(0,C)$ and up to two coexistence equilibrium points $P_1=(u_1,u_1+C)$ and $P_2=(u_2,u_2+C)$, where $u_1\leq u_2$ are given by
\begin{equation}\label{delta}
u_{1,2} = \frac{1}{2}\left(H_1 \pm \sqrt{\Delta}\right)\quad\text{with}\\
\Delta=H_1^2+4H_2,\quad H_1=1-A-Q \quad\text{and}\quad H_2=A-CQ.
\end{equation}
Note that, depending on the system parameters, $u_1$ can be negative and both $u_1$ and $u_2$ can also be complex. In \cite{arancibia,arancibia2,arancibia5,arancibia9} the authors proved that $(0,0)$ is always unstable, $(1,0)$ is always a saddle point and the stability of $(0,C)$ depends on the value of $H_2$ \eref{delta}. If $H_2<0$, then the equilibrium point $(0,C)$ is a saddle point, the equilibrium point is a saddle node if $H_2=0$, and the equilibrium point is a stable node if $H_2>0$. The equilibrium point $P_1$ is always a saddle point (when it is located in the firs quadrant), while $P_2$ can be a stable or unstable node, see \Tref{T01}. In particular, $P_1$ and $(0,C)$ exchange stability by increasing $H_2$ through $\Delta=0$. Moreover, the authors proved that all solutions of \eref{ODEeq03} which are initiated in $\mathbb{R}^2_{\geq0}$ end up in the region
\begin{equation}\label{phi}
\Phi=\{(u,v),\ 0\leq u\leq1,\  0\leq v\leq1+C\}.
\end{equation}
\begin{table}
\begin{tabular}{l l l l l l l}
\hline
$H_1$ & $H_2$ & $\Delta$ & Location of $P_1$ and $P_2$  & $P_2$ is stable if & $P_2$ is unstable if \\
\hline
$>0$	 & $<0$ & $>0$	& $P_1,P_2\in\Phi$ & $S>\frac{\left(H_1+\sqrt{\Delta}\right)\left(Q-\sqrt{\Delta}\right)}{H_1+2A+\sqrt{\Delta}}$ & $S<\frac{\left(H_1+\sqrt{\Delta}\right)\left(Q-\sqrt{\Delta}\right)}{H_1+2A+\sqrt{\Delta}}$\\
$ \neq0$ & $>0$ & $>0$ & $P_1\notin\Phi,P_2\in\Phi$ & &\\ 
\hline
$>0$	 & $<0$ & $=0$ & $P_1=P_2\in\Phi$	& $S>\frac{QH_1}{H_1+2A}$ & $S<\frac{QH_1}{H_1+2A}$\\
\hline
$>0$ & $=0$ & $>0$ & $P_1=(0,C),P_2\in\Phi$ & $S>\frac{H_1\left(Q-H_1\right)}{1-Q}$ & $S<\frac{H_1\left(Q-H_1\right)}{1-Q}$\\
\hline
$=0$ & $>0$ & $>0$ & $P_1\notin\Phi,P_2\in\Phi$ & $S>\frac{\sqrt{H_2}\left(Q-2\sqrt{H_2}\right)}{A+\sqrt{H_2}}$ & $S<\frac{\sqrt{H_2}\left(Q-2\sqrt{H_2}\right)}{A+\sqrt{H_2}}$\\
\hline
$>0$ & $<0$ & $<0$ & \multicolumn{3}{l}{System \eref{ODEeq03} does not have equilibrium points in $\Phi$}\\ 
$\leq0$ & $\leq0$ & & & &\\
\hline
\end{tabular}
\caption{Stability of the coexistence equilibrium point $P_2=(u_2,u_2+C)$ of system \eref{ODEeq03} as derived in \cite{arancibia,arancibia2,arancibia5,arancibia9}, with $H_1$,  $H_2$ and $\Delta$ defined in \eref{delta} and $\Phi$ in \eref{phi}.} \label{T01}
\end{table}

\begin{figure}
\begin{center}
\includegraphics[width=15cm]{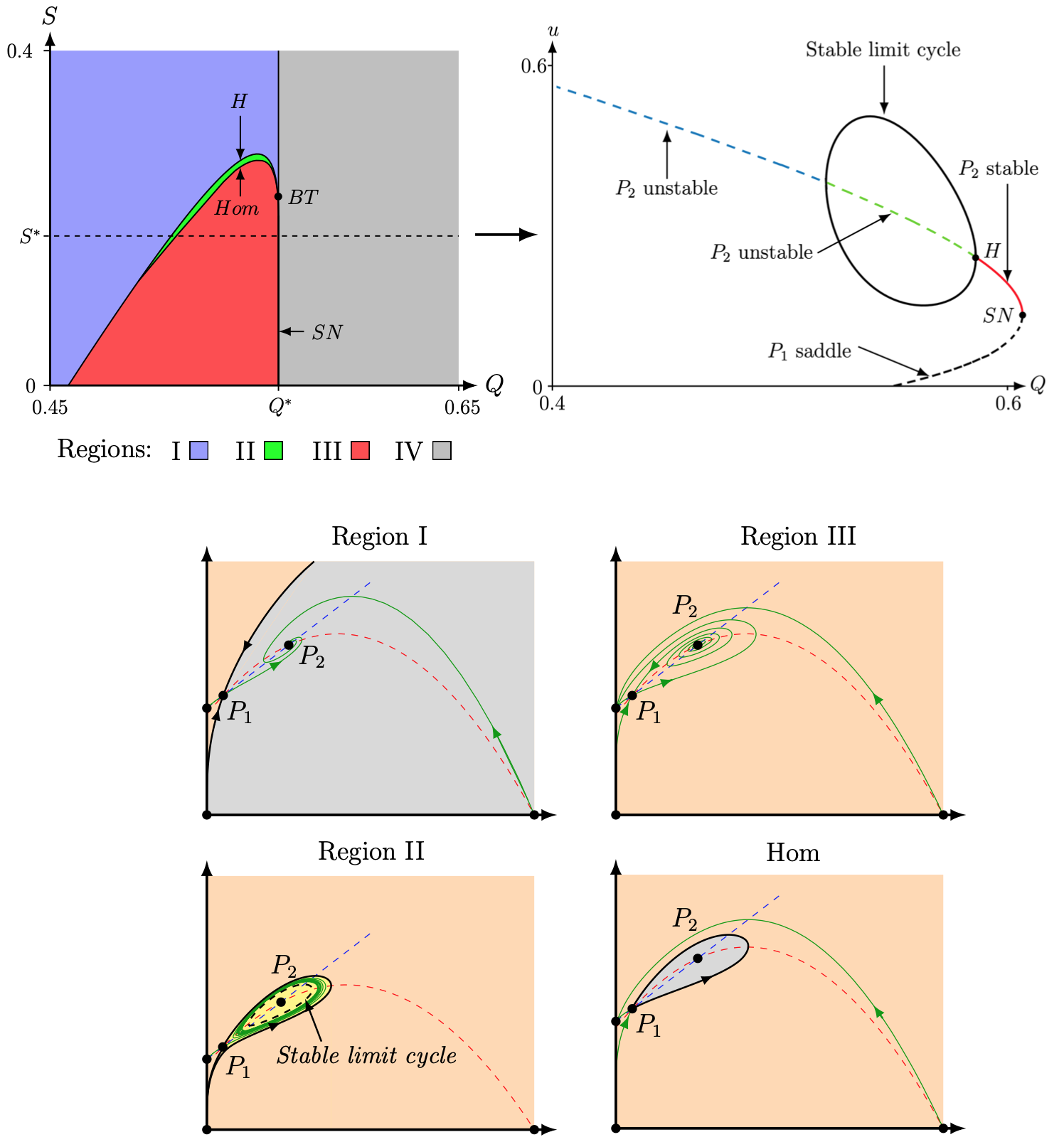}
\end{center}
\caption{In the top left panel we show the bifurcation diagram of system \eref{ODEeq03} for $(A,C)=(0.15,0.28)$ fixed and created with the numerical bifurcation package MATCONT \cite{matcont}. The curve $H$ represents the Hopf curve where $P_2$ changes stability and where a limit cycle is created, $Hom$ represents the homoclinic bifurcation where the limit cycle is destroyed, $SN$ represents the saddle-node curve where $\Delta=0$ and $BT$ represents the Bogdanov--Takens bifurcation. In the top right panel we show the bifurcation diagram of system \eref{ODEeq03} for $(A,C,S^*)=(0.15,0.28,0.23840712)$ fixed and varying the parameter $Q$. We show the behaviour of the equilibrium points $P_1$ and $P_2$ for different values of the parameter $Q$. In the phase plane of system \eref{ODEeq03} (middle and bottom panels) the orange regions represent the basin of attraction of the equilibrium point $(0,C)$, the light blue region represent the basin of attraction of the equilibrium point $P_2$ and the yellow region represent the basin of attraction of a stable limit cycle (only presented in Region II).}
\label{F01}
\end{figure}

From \eref{delta} and \Tref{T01} we can conclude that a modification of the parameter $Q$ changes the location of the equilibrium points $P_1$ and $P_2$ and this variation also changes the stability of the equilibrium point $(0,C)$ and $P_2$. Moreover, the variation of the parameter $S$ changes the stability of the equilibrium point $P_2$. Therefore, the basins of attraction of the equilibrium points $(0,C)$ and $P_2$ depend on the parameters $Q$ and $S$. In Figure \ref{F01} the numerical bifurcation package MATCONT \cite{matcont} is used to obtain the bifurcation diagram\footnote{The Matlab package ode45 was used to generate the data for the simulations and then the PGF package (or tikz) was used to generate the graphics format.} of the diffusion-free modified Holling--Taner model  \eref{ODEeq03} for the system parameters $(A,C)=(0.15,0.28)$ fixed, see top panel of \Fref{F01}. We choose the parameter values $(A,C)=(0.15,0.28)$ fixed since these are the values used in \cite{arancibia,arancibia2,arancibia5,arancibia9}. Alternatively, we could have fixed $Q$ and changed $A$ or $C$ and this would have resulted in equivalent bifurcation diagrams. The bifurcation curves divide the $(Q,S)$-space into four different areas with different behaviour: 
\begin{itemize}
\item \textbf{Region I ($\Delta>0$)}: the equilibrium points $(0,C)$ and $P_2$ are stable nodes and the stable manifold of $P_1$ determines the boundary of the domain of attraction of $(0,C)$ and the domain of attraction of $P_2$ (orange and light blue regions in \Fref{F01} respectively). All initial conditions initiated above this separatrix go to $(0,C)$, which represents the scaled alternative food level and represents the extinction of the prey and the persistence of the predator population, while all solutions which are initiated below the separatrix go to $P_2$ which represent the stabilisation of both populations. 
\item \textbf{Region II ($\Delta>0$)}: the equilibrium point $(0,C)$ is still a stable node, while $P_2$ is now an unstable node surrounded by a stable limit cycle. This limit cycle is born through a Hopf bifurcation and terminated by a homoclinic bifurcation. The stable manifold of $P_1$ again determines the boundary of the domain of attraction of $(0,C)$ and the domain of attraction of the stable limit cycle (yellow region in \Fref{F01}). The stable limit cycle represents the oscillation of both populations.
\item \textbf{Region III ($\Delta>0$)}: the limit cycle is terminated and the equilibrium point $(0,C)$ is a global stable node and $P_2$ is an unstable node. Therefore, all initial conditions go to $(0,C)$ and hence the prey goes extinct.
\item \textbf{Region IV ($\Delta<0$)}: system \eref{ODEeq03} does not have equilibrium points in the first quadrant and $(0,C)$ is a global stable node. Similarly with \textbf{Region III}, all trajectories with initial conditions in the first quadrant go to $(0,C)$ and hence the prey goes extinct.
\item \textbf{$Q=Q^*$ ($\Delta=0$)}: this form the boundary between \textbf{Regions I, II \& III} and \textbf{Region IV}. On this line the equilibrium point $(0,C)$ is a stable node and the equilibrium points $P_1$ and $P_2$ collapse. So, system \eref{ODEeq03} experiences a saddle-node bifurcation (labeled $SN$ in Figure \ref{F01}) and a Bogdanov--Takens bifurcation (labeled $BT$ in \Fref{F01}) along this line when $S=QH_1/(H_1+2A)$ \cite{arancibia9}. 
\end{itemize}

\section{Spatio-temporal Model}\label{stem}

In this section, we present the model with diffusion, we first recall the criteria for Turing instability for a general spatio-temporal predator-prey model (where the populations are assumed to be distributed in an unbounded domain)
\begin{equation}\label{GENPDE}
\begin{aligned}
\frac{\partial A}{\partial t} &=  W\left(A,B\right)+\triangledown^2A,\\
\frac{\partial B}{\partial t} &= Z\left(A,B\right)+d\triangledown^2B.
\end{aligned}
\end{equation}
Here, $A(X,t)$ and $B(X, t)$ are considered to be the prey and the predator population respectively, $W(A,B)$ and $Z(A,B)$ describe their nonlinear interaction and $d=D_A/D_B$ with $D_A$ and $D_B$ constant diffusivities. Turing \cite{turing} showed that an equilibrium point that is stable in a temporal model can become unstable upon adding diffusion in the model. In the absence of diffusion, we analyse the stability of an equilibrium point $\left(A_0,B_0\right)$ such that $W\left(A_0,B_0\right)=Z\left(A_0,B_0\right)=0$. The stability of this equilibrium point depends on the eigenvalues of the Jacobian matrix $J(A_0,B_0)$ which can be found by solving $|J-\lambda I|=0$. That is, $\lambda^2-\left(W_A+Z_B\right)\lambda+\left(W_AZ_B-W_BZ_A\right)=0$, where $W_A=\partial W/\partial A$, $W_B=\partial W/\partial B$, $Z_A=\partial Z/\partial A$ and $Z_B=\partial Z/\partial B$ are evaluated at the equilibrium point $(A_0,B_0)$. The equilibrium point $\left(A_0,B_0\right)$ is stable if
\begin{equation}\label{ODE}
\begin{aligned}
\tr(J(A_0,B_0)) &=W_A+Z_B<0\quad\text{and}\quad\det(J(A_0,B_0)) &=W_AZ_B-W_BZ_A>0.
\end{aligned}
\end{equation}
The Turing instability is obtained by linearising the PDE system around the equilibrium point $\left(A_0,B_0\right)$. The stability of the equilibrium point is now determinated by the roots of the characteristic polynomial $|J-\lambda I-k^2D|$ where $k$ is the wave number \cite{malchow} and 
$D=\left(\begin{array}{cc}1 & 0 \\ 0 & d\end{array}\right)$. 
This defines the dispersion relation $\lambda\left(k\right)$ which is the solution of $\lambda^2-\alpha\left(k^2\right)\lambda+\beta\left(k^2\right)=0$ where $\alpha\left(k^2\right)=\tr(J\left(A_0,B_0\right)-k^2\left(1+d\right)$ and $\beta\left(k^2\right)=dk^4-\left(dW_A+Z_B\right)k^2+\det(J\left(A_0,B_0\right))$. Hence, we obtain the dispersion relation $\lambda(k)$ 
\begin{equation}\label{lambda}
\lambda_{\pm}(k)=\frac{1}{2}\left(\alpha\left(k^2\right)\pm\sqrt{\left(\alpha\left(k^2\right)\right)^2-4\beta\left(k^2\right)}\right).
\end{equation}
If we assume that the conditions defined in \eref{ODE} hold then $\alpha\left(k^2\right)<0$ and the equilibrium point $\left(A_0,B_0\right)$ is thus unstable in \eref{GENPDE} if we have that $\beta\left(k^2\right)<0$. By \eref{ODE} we also have that $\det(J(A,B))>0$ and the minimum of the quadratic $\beta(k^2)$ occurs when $k^2=\left(dW_A+Z_B\right)/\left(2d\right)$. Therefore, the conditions for the equilibrium point $\left(A_0,B_0\right)$ to be unstable in \eref{GENPDE} are
\begin{equation}\label{PDE}
dW_A+Z_B>0\quad\text{and}\quad \left(dW_A+Z_B\right)^2-4d\det(J\left(A_0,B_0\right))>0.
\end{equation}
The Turing conditions in two-dimensional space can be obtained from \eref{ODE} and \eref{PDE} by replacing $k^2$ with $m^2+\ell^2$, where $(m,\ell)\in\mathbb{R}^2$ indicate the wave numbers in the $x$-direction and $y$-direction respectively \cite{banerjee1}.

\subsection{Equilibrium point $P_2$}\label{nusip2}
We now discuss the Turing conditions \eref{ODE} and \eref{PDE} for the model of interest in this manuscript and we first focus on the coexistence equilibrium point $P_2$. The only other equilibrium point that can be stable in the temporal system is $(0,C)$, we will investigate its Turing space in subsection \ref{epC}. These Turing conditions refer to the equilibrium point $P_2$ to be stable in the diffusion-free system \eref{ODEeq03} and unstable in the full system \eref{eq03}. The conditions \eref{ODE} for the equilibrium point $P_2$ to be stable in the diffusion-free system are met if we assume that $H_1>0$, $H_2<0$ and $\Delta>0$ (see \Tref{T01}). That is,   
\begin{equation}\label{p2ODE1}
\frac{\left(H_1+\sqrt{\Delta}\right)\left(Q-\sqrt{\Delta}\right)}{\left(H_1+2A+\sqrt{\Delta}\right)}-S<0~\text{and}~\frac{S\left(H_1+\sqrt{\Delta}\right)\sqrt{\Delta}}{\left(H_1+2A+\sqrt{\Delta}\right)}>0.
\end{equation}
The conditions \eref{PDE} for the equilibrium point $P_2$ to be unstable in \eref{eq03} are
\begin{equation}\label{p2PDE1}
\begin{aligned}
&d\frac{\left(H_1+\sqrt{\Delta}\right)\left(Q-\sqrt{\Delta}\right)}{\left(H_1+2A+\sqrt{\Delta}\right)}-S>0~\text{and}\\
&\left(\frac{d\left(H_1+\sqrt{\Delta}\right)\left(Q-\sqrt{\Delta}\right)}{\left(H_1+2A+\sqrt{\Delta}\right)}-S\right)^2-\frac{4dS\left(H_1+\sqrt{\Delta}\right)\sqrt{\Delta}}{\left(H_1+2A+\sqrt{\Delta}\right)}>0.
\end{aligned}
\end{equation}
Note that these assumptions on $H_1$, $H_2$ and $\Delta$  are the most general case shown in \Tref{T01} for which conditions \eref{ODE} are met. That being said, there are other cases for which these conditions are also met, for brevity, we will not investigate these cases. Additionally, the other coexistence equilibrium point $P_1$ never fulfils \eref{ODE} as it is a saddle point in the temporal system.

\begin{figure}
\begin{center}
\includegraphics[width=15.5cm]{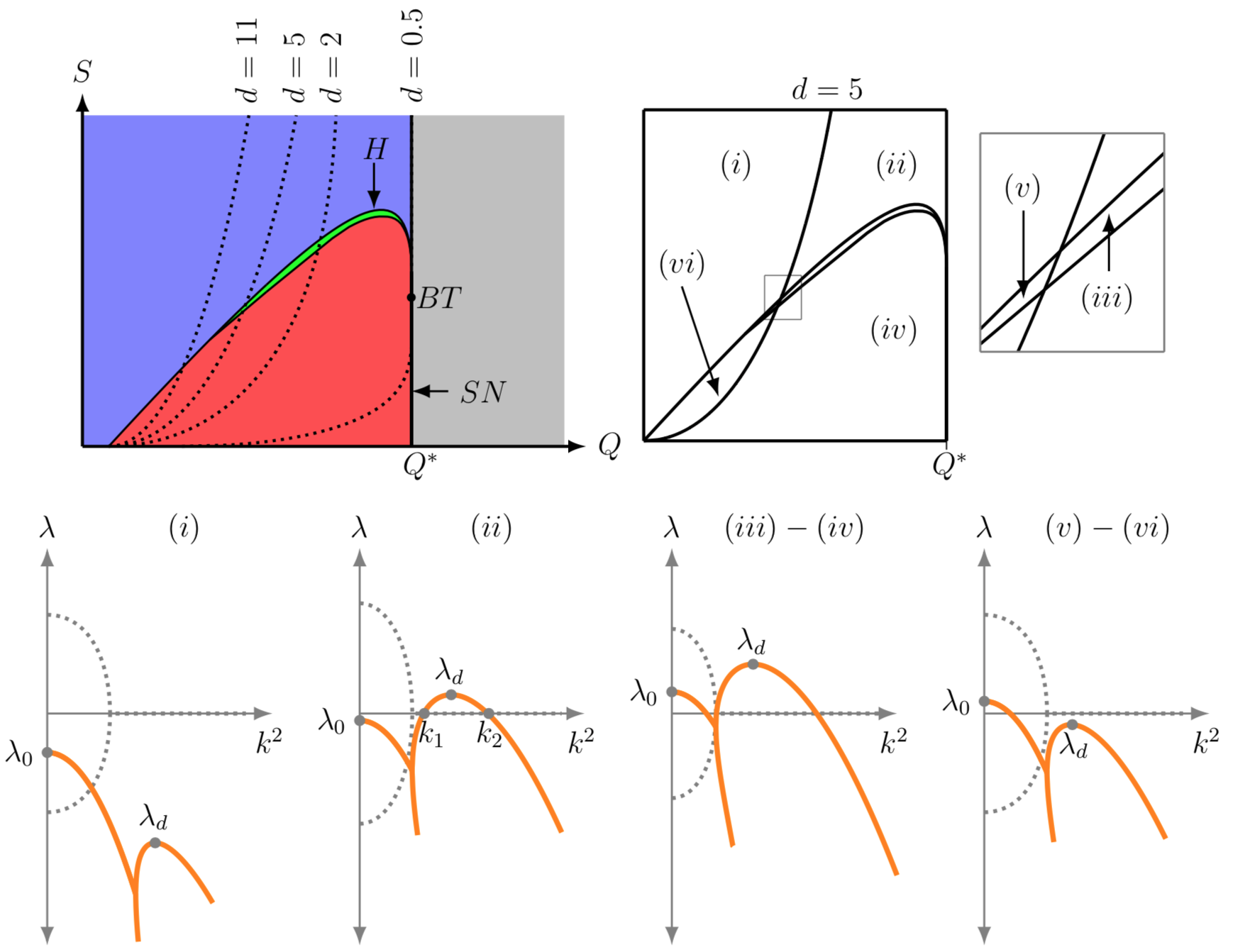}
\end{center}
\caption{In the top panel we show the Turing space of system \eref{eq03} for $P_2$ with $(A,C)=(0.15,0.28)$ fixed. In the left top panel we consider four different values of the diffusion ratio $d$ $(0.5,2,5,11)$ (dotted curves) and in the right top panel we consider the case where the diffusion ratio $d=5$. In the bottom panel, we show the real part (orange curve) and imaginary part (grey dotted curve) of the dispersion relation $\lambda(k)$ \eref{lambda} as a function of the wave number squared for the system parameter $(A,C,d)=(0.15,0.28,5)$ fixed and (i) if $Q=0.5$ and $S=0.27$ then $P_2$ is stable in the ODE and PDE system; (ii) if $Q=0.575$ and $S=0.26$ then $P_2$ is stable in the ODE and unstable in the PDE system; (iii) if $Q=0.575$ and $S=0.235$ or (v) if $Q=0.52$ and $S=0.112$ then $P_2$ is unstable surrounded by a stable limit cycle in the ODE system and unstable in the PDE system; (iv) if $Q=0.575$ and $S=0.1$ or (vi) if $Q=0.51$ and $S=0.18$ then $P_2$ is unstable without a limit cycle in the ODE system and unstable in the PDE system.}
\label{F02}
\end{figure}

As before, we fix the system parameter $(A,C)=(0.15,0.28)$ in system \eref{eq03} and the Turing parameter space $(Q,S)$ of the equilibrium point $P_2$ is given in \Fref{F02}. Here, the conditions for the diffusion-driven instability \eref{p2ODE1} and \eref{p2PDE1} are met in region $(ii)$ of \Fref{F02} and we thus expect to see Turing patterns in this region. In \Tref{T02} we summarise the stability properties of the equilibrium point $P_2$ in system \eref{eq03} with and without diffusion.
\begin{table}
\begin{tabular}{l l l l}
\hline
Region 	& $\lambda_0$ 	& $\lambda_d$ 	& \\
\hline
$(i)$ 		& $<0$ 	& $<0$ 	& $P_1\in\Phi$ and $P_2$ is stable in the ODE and the PDE system  				\\
\hline
$(ii)$ 	& $<0$ 	& $>0$  	& $P_1\in\Phi$ and $P_2$ is stable in the ODE system and unstable 		\\
 		&  		&   		& in the PDE system (Turing patterns) 												\\
\hline
$(iii)$ 	& $>0$ 	& $>0$ 	& $P_1\in\Phi$ and $P_2$ is unstable surrounded by a stable limit cycle 		\\
 		&  		&  		& in the ODE system and unstable in the PDE system										\\
\hline
$(iv)$ 	& $>0$ 	& $>0$ 	& $P_1\in\Phi$ and $P_2$ is unstable without limit cycle in the ODE					\\
 		&  		&  		&  system and unstable in the PDE system 									\\
\hline
$(v)$ 	& $>0$ 	& $<0$ 	& $P_1\notin\Phi$ and $P_2$ is unstable surrounded by a stable limit cycle	\\
 		&  		&  		& in the ODE system and unstable in the PDE system 									\\
\hline
$(vi)$ 	& $>0$ 	& $<0$ 	& $P_1\notin\Phi$ and $P_2$ is unstable without limit cycle in the ODE  				\\
 		&  		&  		& system and unstable in the PDE system 									\\
\hline
\end{tabular}
\caption{Summary of the dispersion relation $\lambda(k)$ \eref{lambda} of system \eref{eq03} for $P_2$ showed in \Fref{F02} for system parameters $(A,C,d)=(0.15,0.28,5)$ fixed. There are six different cases, depending on $\lambda_0=\mathbf{Re}(\lambda(0))$ and $\lambda_d=max(\mathbf{Re}(\lambda(k)))$ and the potential limit cycle in the temporal system.}\label{T02}
\end{table} 
\subsubsection{Numerical Simulations near the equilibrium point $P_2$}\label{2Dp2}
Here, we present numerical solution of system \eref{eq03} in one and two-dimensional space for the system parameters $(A,C)=(0.15,0.28)$ fixed and for initial conditions near $P_2=(u_2+u_2+C)$. In particular, the initial condition is 
\begin{equation}\label{ic}
u_0=u_2+0.012e^{-7x^2}\quad\text{and}\quad v_0=u_0+C. 
\end{equation} 
The numerical integration of system \eref{eq03} is performed by using the Matlab PDEPE toolbox with zero-flux boundary conditions on a domain of size $200$ and we discretise using $1500$ grid points. Note that increasing the domain size upon which the simulations were performed did not significantly change the observed results.

In region $(i)$, see \Fref{F02} and \Tref{T02}, the equilibrium point $P_2$ is stable in the temporal system and in the spatio-temporal system. We observe that the initial condition \eref{ic} evolves as expected to the spatially homogeneous stationary equilibrium point $P_2$ (the simulation is not shown). 

\begin{figure}
\begin{center}
\includegraphics[width=16cm]{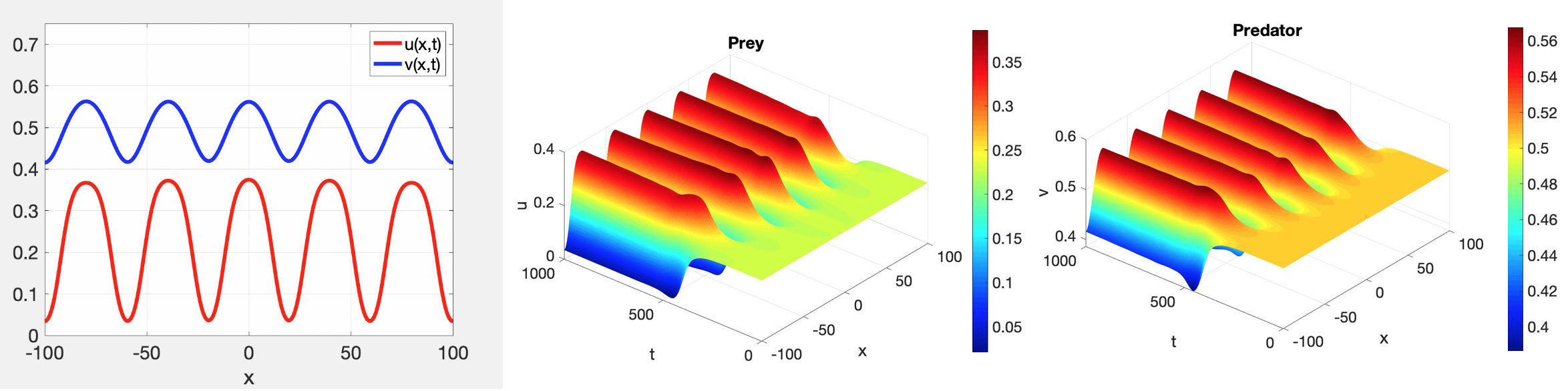}
\end{center}
\caption{Numerical simulation of system \eref{eq03} in one-dimensional space with system parameters $(A,C,d,Q,S)=(0.15,0.28,5,0.575,0.26)$ and initial condition defined in \eref{ic}. For the parameter values the equilibrium point $P_2=(0.22642,0.50642)$ is unstable in the spatio-temporal system but stable in the temporal system and the formation of Turing patterns is expected, see $(ii)$ in \Fref{F02}. We observe a Turing pattern that is stationary in time and oscillatory in space. An animated version of this figure is accessible on \href{http://www.doi.org/10.6084/m9.figshare.10059242}{http://www.doi.org/10.6084/m9.figshare.10059242}.}
\label{F03}
\end{figure}

\begin{figure}
\begin{center}
\includegraphics[width=7.5cm]{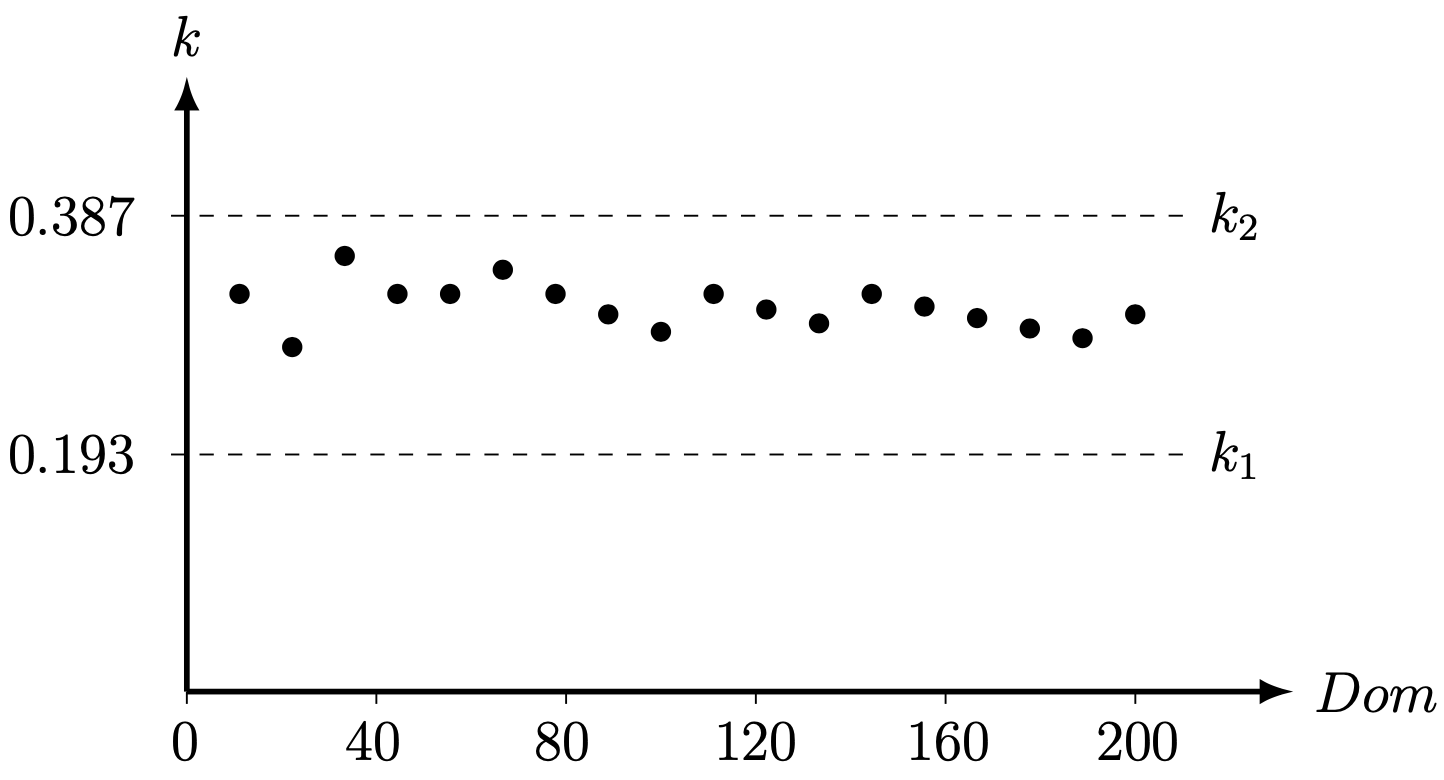}
\end{center}
\caption{The numerical observed wave number of the Turing pattern as function of the domain size ($Dom$) for system parameters $(A,C,d,Q,S)=(0.15,0.28,5,0.575,0.26)$ fixed. We see that the observed wave numbers are in between the minimum ($k_{1}$) and maximum ($k_{2}$) for which the stable homogeneous solutions becomes unstable when diffusion is included in the model, see \Fref{F02}.}
\label{F04}
\end{figure}

In region $(ii)$, see \Fref{F02} and \Tref{T02}, the equilibrium point $P_2$ is stable in the temporal system and unstable in spatio-temporal system and the conditions for a Turing instability are thus met. Indeed, the initial condition \eref{ic} evolves to a Turing pattern that is a periodic solution in space around $P_2$ and that is stationary in time, see \Fref{F03}. Note that the simulation is run over a long period of time to ensure that the Turing pattern is stationary in time. Note that the analysis of the Turing instability of \Sref{nusip2} was based on an unbounded domain, while the numerical simulation are performed on a bounded domain. The analysis on the unbounded domain shows that there is a range of unstable wave numbers, i.e. $\mathbf{Re}(\lambda(k))>0$ for $k\in(k_{1},k_{2})$ in region $(ii)$, see \Fref{F02}. On a bounded domain a similar Turing analysis can be done by taking the solution as $u(x,t)=\alpha e^{\lambda t}cos(kx)$ where $\alpha$ is the initial amplitude and $k$ is the wave number. This results in a discrete set of (unstable) wave numbers lying on $\lambda(k)$ and depending on the domain size, i.e. the spatial period has to fit in the domain. In \Fref{F04} we show the numerical observed wave number of the Turing pattern as function of the domain size and we see that the observed wave number is as expected in between $k\in(k_{1},k_{2})$.    

\begin{figure}
\begin{center}
\includegraphics[width=16cm]{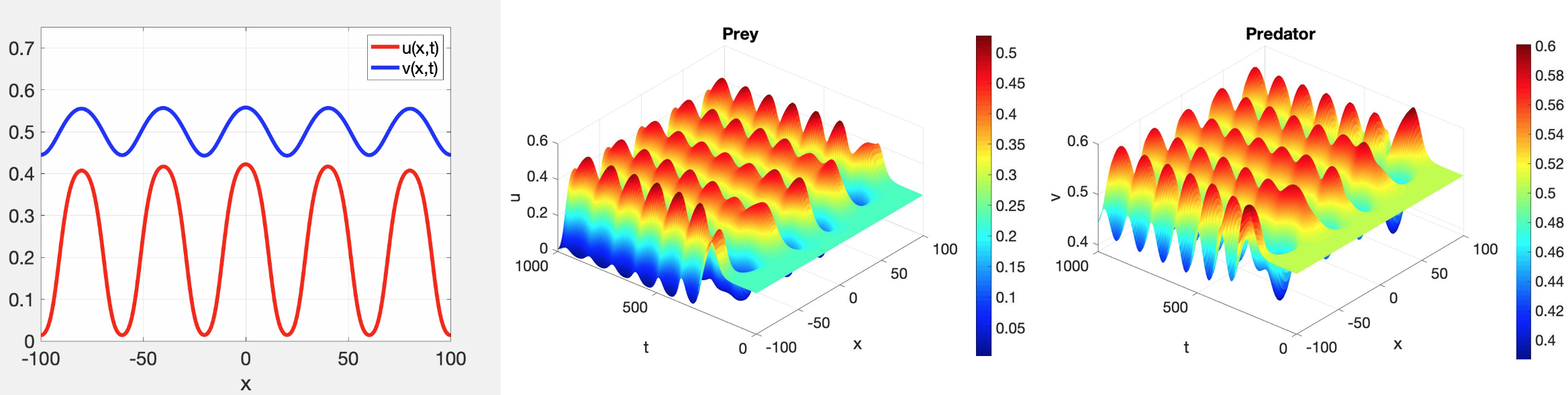}
\includegraphics[width=16cm]{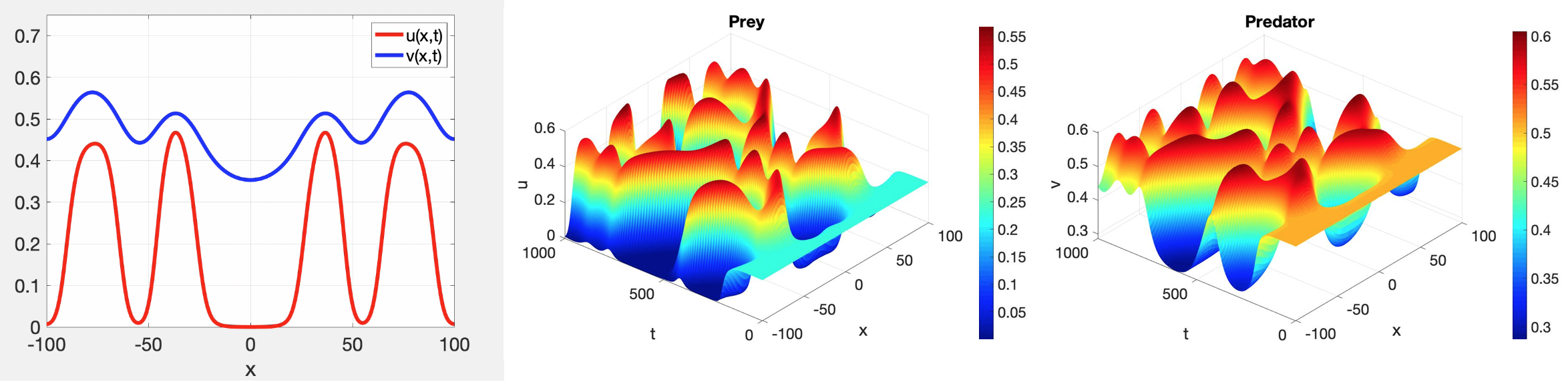}
\end{center}
\caption{Numerical simulation of system \eref{eq03} in one-dimensional space with system parameters $(A,C,d)=(0.15,0.28,5)$ fixed and initial condition defined in \eref{ic}. In the top panel $(Q,S)=(0.575,0.1)$ and the equilibrium point $P_2=(0.22642,0.50642)$ is unstable surrounded by stable limit cycle in the ODE system and unstable in the PDE system, see $(iii)$ in \Fref{F02}. We observe that the solution is oscillatory in space and in time with period $113.63(t)$ and $40(x)$ respectively. An animated version of this figure is accessible on \href{http://www.doi.org/10.6084/m9.figshare.10059248}{http://www.doi.org/10.6084/m9.figshare.10059248}. In the bottom panel $(Q,S)=(0.575,0.07)$ and the equilibrium point $P_2$ is unstable and there is no limit cycle in the ODE, see $(iv)$ in \Fref{F02}. We observe that the solution evolves to a different pattern observed in the top panel that is less regular pattern. An animated version of this figure is accessible on \href{http://www.doi.org/10.6084/m9.figshare.10059254}{http://www.doi.org/10.6084/m9.figshare.1005925}.} 
\label{F05}
\end{figure}

In region $(iii)$ and $(iv)$, see \Fref{F02} and \Tref{T02}, the equilibrium point $P_2$ is unstable with respect to wave numbers near to zero and near to $\lambda_d$. In region $(iii)$ the equilibrium point $P_2$ is  surrounded by a stable limit cycle in the temporal system, while this is not the case in region $(iv)$.  Note that in these two regions the equilibrium point $P_1$ is located in $\Phi$ and thus $(0,C)$ is a stable node in the temporal system. In region $(iii)$, we observe that the initial condition \eref{ic} evolves to pattern that is oscillatory in space and in time, see top panel of \Fref{F05}. Additionally, we observe that in \Fref{F05} the period of the oscillation in the spatio-temporal system is $40(x)$ and the wave number is $k=0.15708$ which is in between the minimum ($k_{1}=0.135$) and maximum ($k_{2}=0.437$) wave number for region $(iii)$, see bottom panel of \Fref{F02}. In region $(iv)$ it is   observed that the initial condition \eref{ic} evolves to a different pattern that is less organised, see bottom panel of \Fref{F05}.

\begin{figure}
\begin{center}
\includegraphics[width=7.5cm]{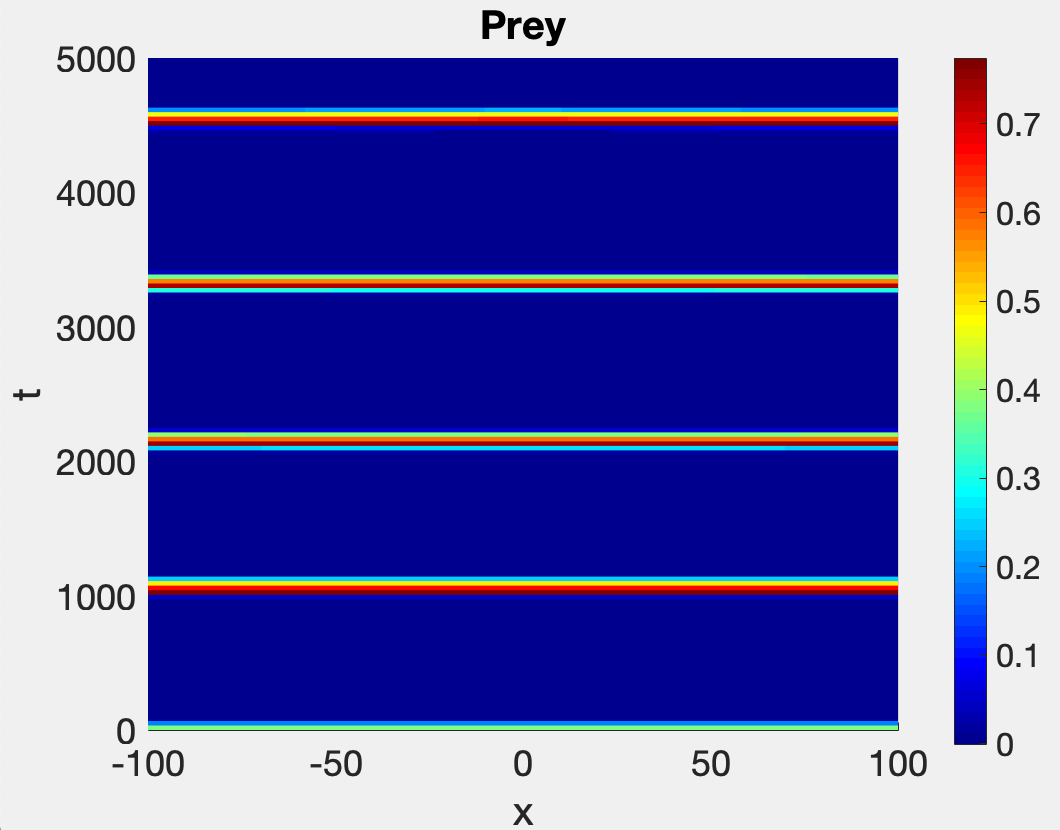}
\includegraphics[width=7.5cm]{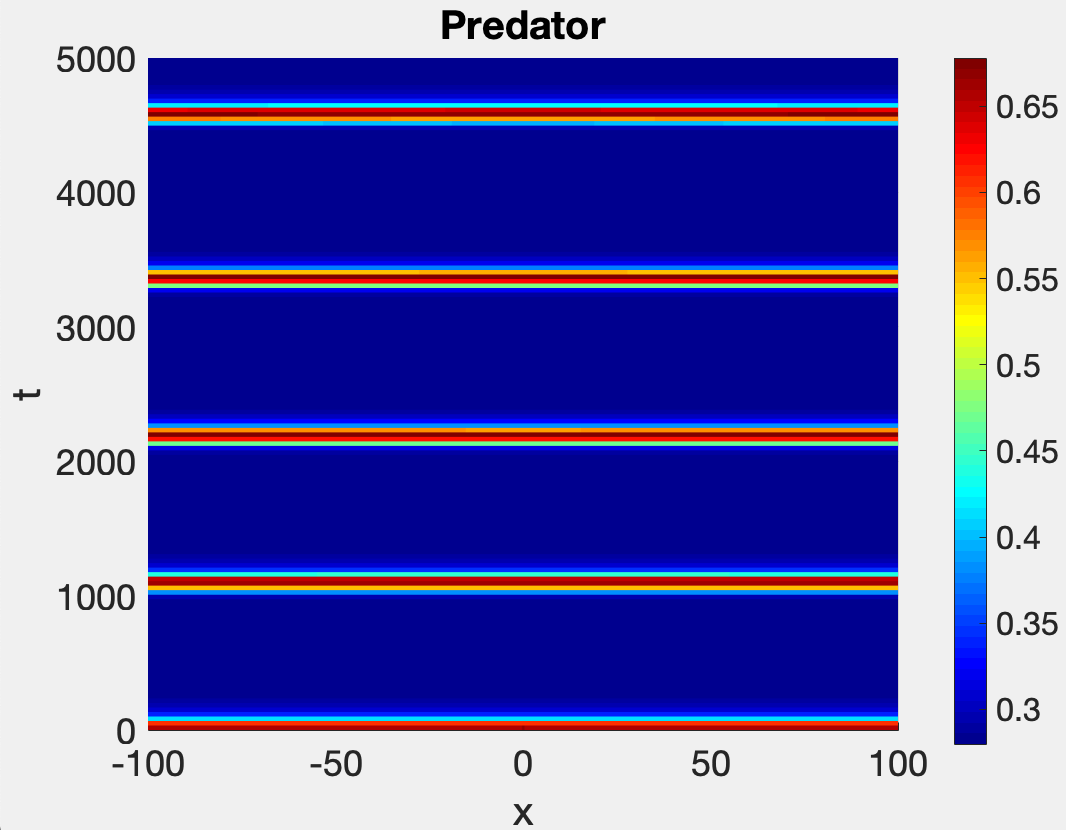}
\includegraphics[width=7.5cm]{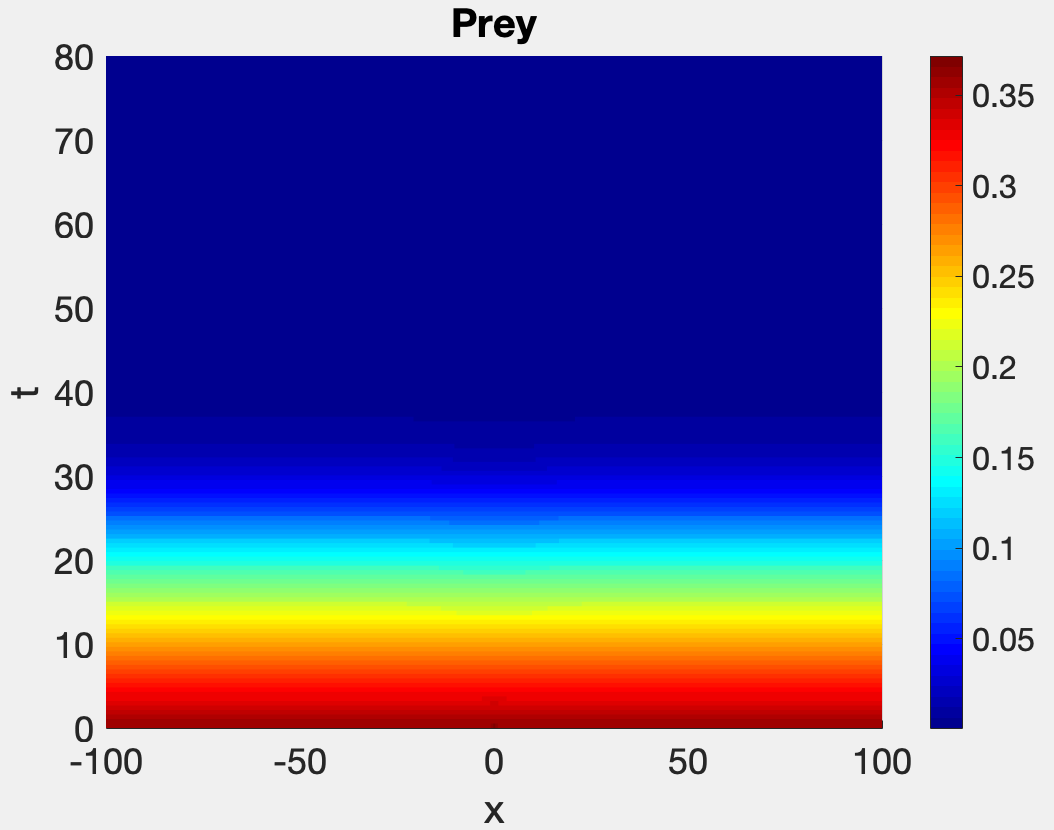}
\includegraphics[width=7.5cm]{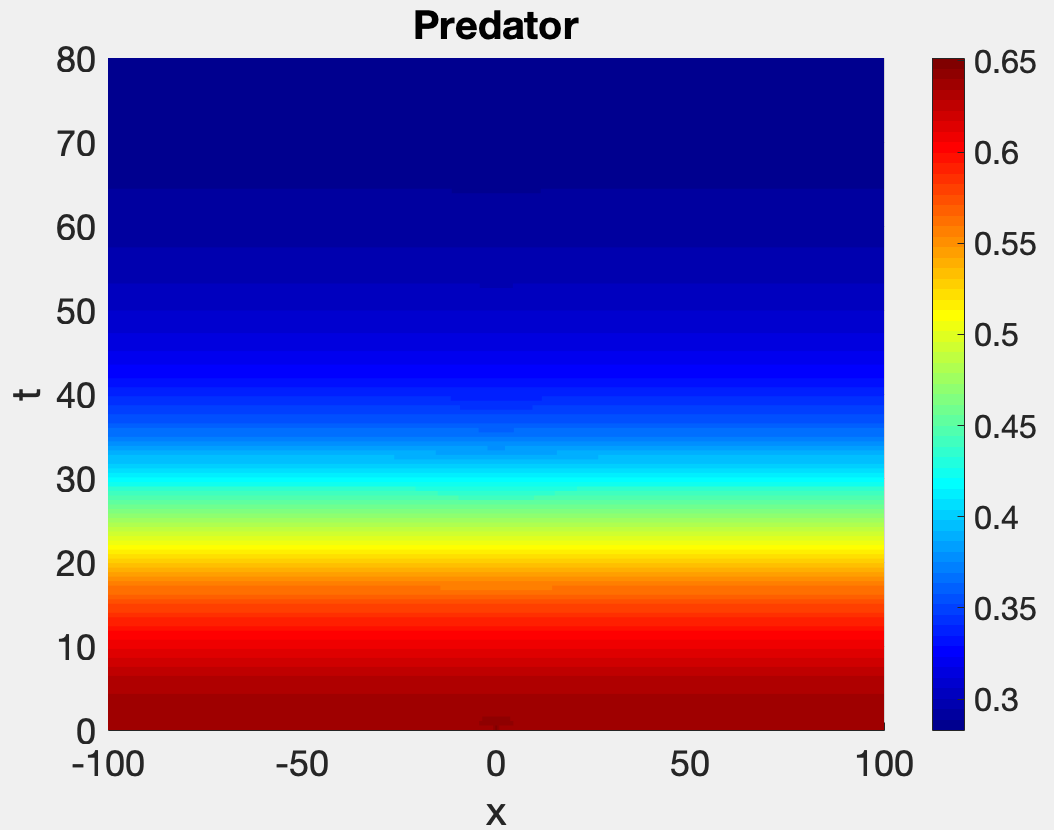}
\end{center}
\caption{Numerical simulation of system \eref{eq03} in one-dimensional space with system parameters $(A,C,d)=(0.15,0.28,5)$ fixed and initial condition defined in \eref{ic}. In the top panel $(Q,S)=(0.52,0.05)$ and the equilibrium point $P_2$ is unstable surrounded by stable limit cycle in the ODE system, see $(v)$ in Figure \ref{F02}. We observe the spatio-temporal pattern that is stationary in space and oscillatory in time as expected since the equilibrium point $P_2$ is only unstable with respect to wave numbers near zero. In the bottom panel $(Q,S)=(0.53,0.085)$ and the equilibrium point $P_2=(0.36,0.64)$ is no longer surrounded by a limit cycle in the ODE system. We observe that the initial condition evolves to $(0,C)$.}
\label{F06}
\end{figure}

In region $(v)$ and $(vi)$, the equilibrium point $P_2$ is only unstable with respect to wave numbers near zero. Additionally, in region $(v)$ the equilibrium point $P_2$ is surrounded by a stable limit cycle in the temporal system, while there is no limit cycle in the temporal system in region $(vi)$. In region $(v)$, we observe that the initial condition evolves to a spatio-temporal pattern that is stationary in space and oscillatory in time, see top panel of \Fref{F06}. We found that the period of the associated limit cycle in the temporal system  is $433.3(t)$ which is close to the period of $463.787(t)$ observed in the top panel of \Fref{F06}. In region $(vi)$, we observe that the solution goes to the equilibrium point $(0,C)$.

\begin{figure}
\begin{center}
\includegraphics[width=7.5cm]{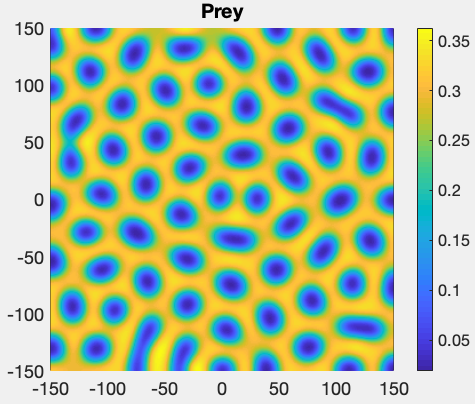}
\includegraphics[width=7.5cm]{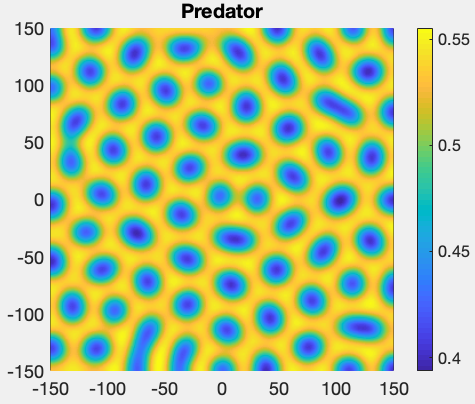}\\
\includegraphics[width=7.5cm]{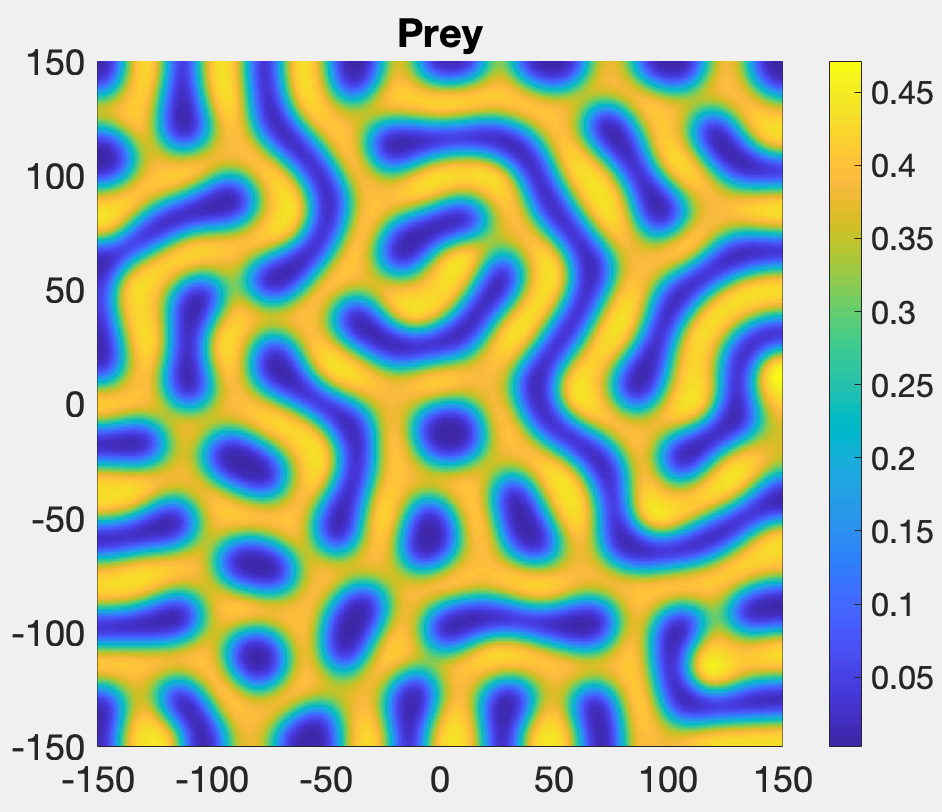}
\includegraphics[width=7.5cm]{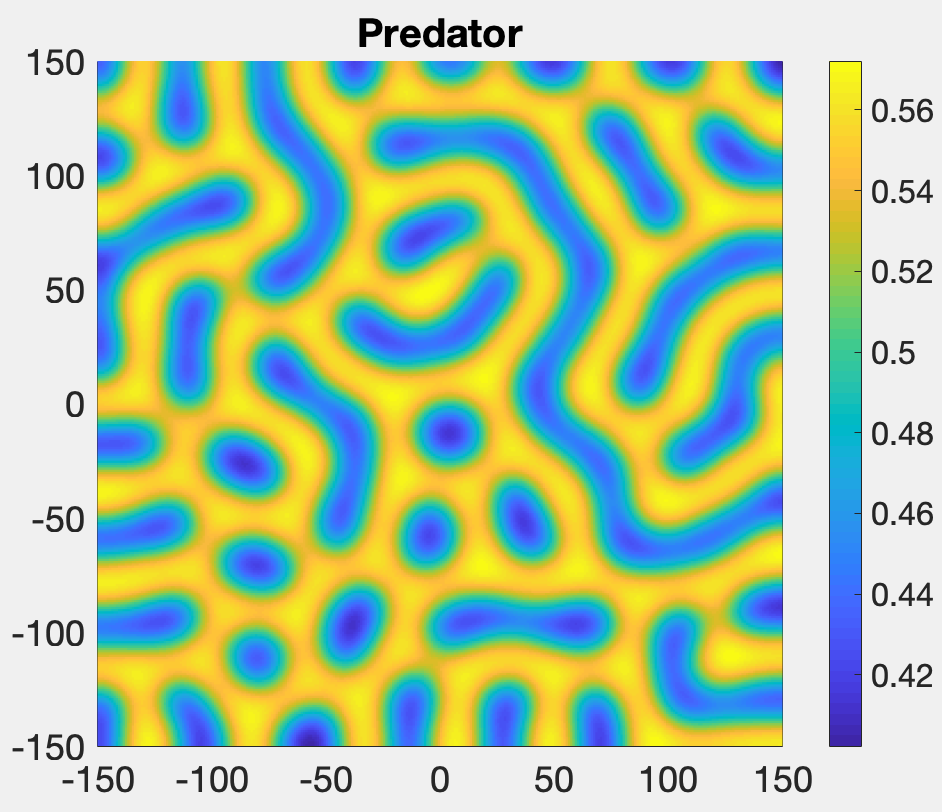}\\
\includegraphics[width=7.5cm]{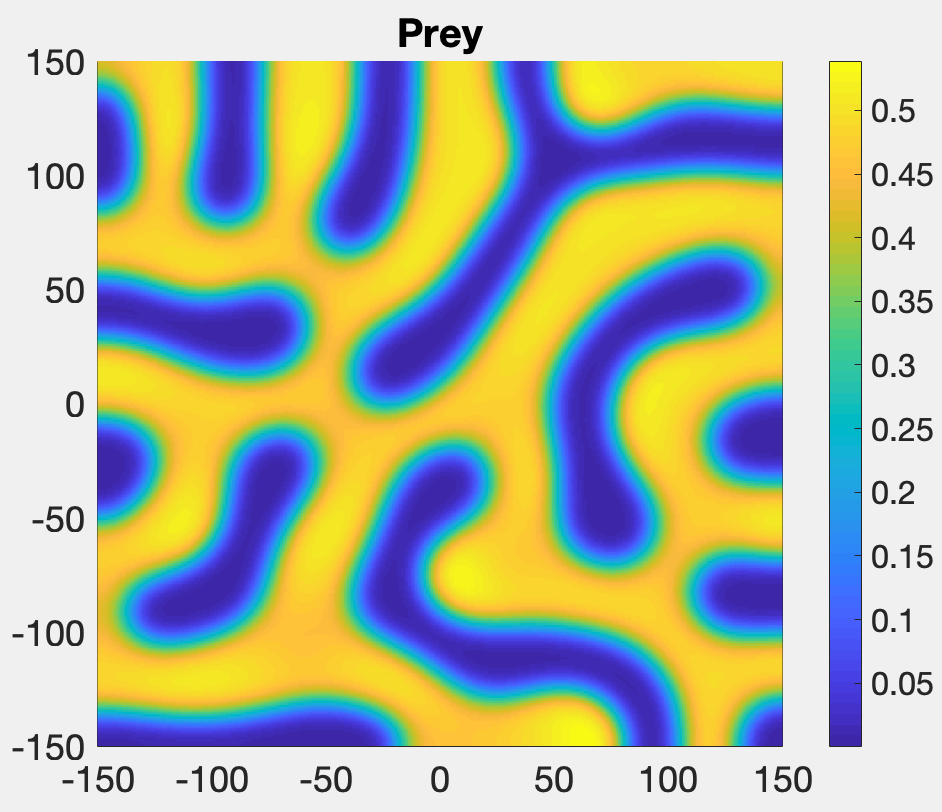}
\includegraphics[width=7.5cm]{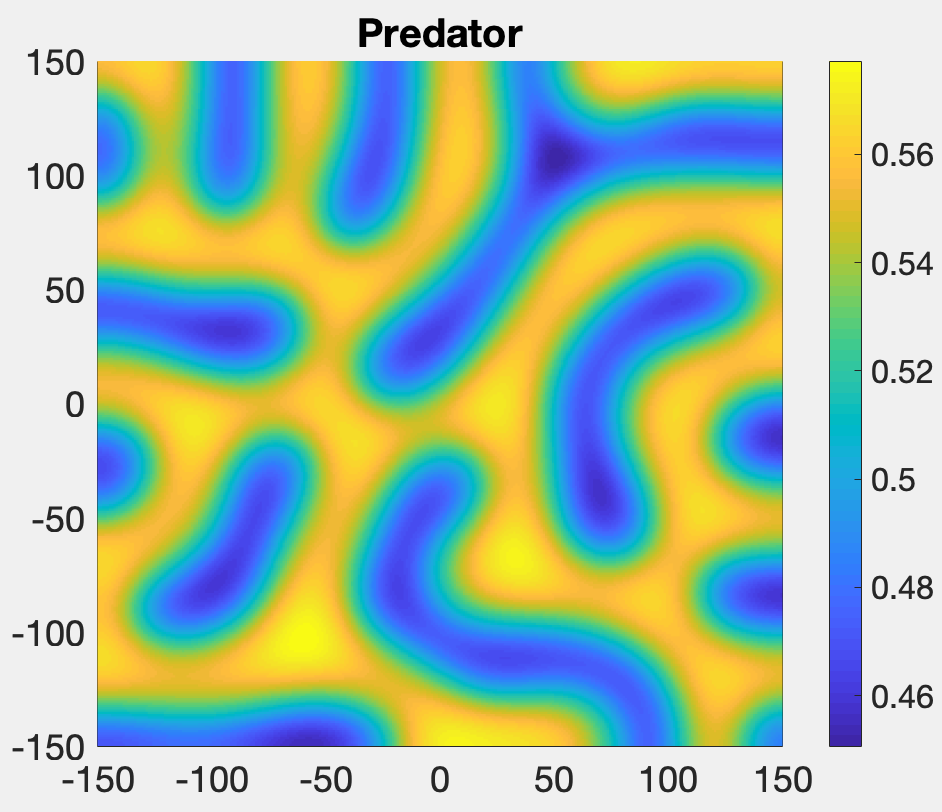}
\end{center}
\caption{Numerical simulation of the predator and prey population for system \eref{eq03} in two-dimensional space for initial conditions near to the equilibrium point $P_2$. We consider system parameters $(A,C,Q,S)=(0.15,0.28,0.575,0.26)$ fixed and $d=5$ in the top panel; $d=9$ in the middle panel; and $d=55$ in the bottom panel. We observe that by increasing the ratio of diffusivities the cold-spot pattern started to coalesce creating combination of labyrinthine and cold-spot pattern. Note that these Turing patterns are stationary in time.}
\label{F07}
\end{figure}

We now shortly discuss the spatio-temporal Turing patterns of system \eref{eq03} in two-dimensional space for region $(ii)$ in \Fref{F02} . We present different Turing patterns by taking system parameter $(A,C,Q,S)=(0.15,0.28,0.575,0.26)$ fixed and varying the ratio of diffusivities. The numerical integration of system \eref{eq03} is performed by using an Euler method for the time integration \cite{garvie,mathews} with a time step size $\Delta t=0.2$ and a finite difference algorithm for a predator-prey system with spatial variation in two-dimensional Laplacian with the zero-flux boundary conditions. The initial condition is a random perturbation around the positive equilibrium point $P_2=(0.22642,0.50642)$. Simulations are run for a long time to ensure that the resulting patterns are stationary in time. We observe that the Turing patterns of the predator and the prey population have the same characteristics. For the ratio of diffusivity $d=5$ we observe a stationary cold-spot pattern over the whole domain, see the top panels in \Fref{F07}. When the ratio of diffusivity is being increased to $d=9$ we observe that the cold-spots started to coalesce creating combination of labyrinthine and cold-spot pattern which coexist in the space, see the middle panels in \Fref{F07}. Finally, by increasing the ratio of diffusivity up to $d=55$ we observe a labyrinthine pattern over the whole domain, see the bottom panel in \Fref{F07}. Moreover, we observe in \Fref{F07} the minimum of the prey population is $0$ while the maximum increases from $0.3$ up to $0.5$ approximately by increasing the ratio of diffusivity. In contrast, the minimum of the predator population increases from $0.43$ up to $0.48$ approximately by also increasing the ratio of diffusivity while the maximum remain constant in approximately $0.55$.
 
\subsection{Equilibrium point $(0,C)$}\label{epC}
\begin{figure}
\begin{center}
\includegraphics[width=9cm]{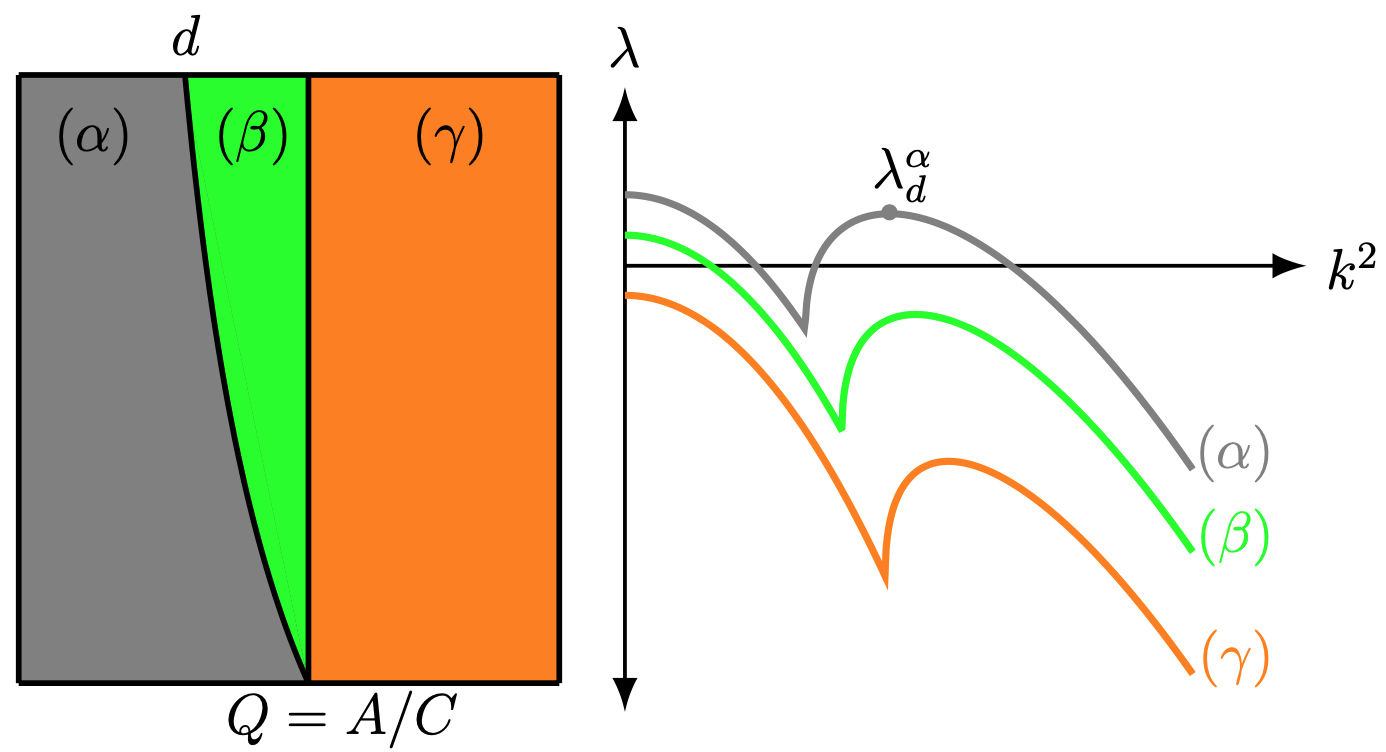}
\end{center}
\caption{In the left panel we show the bifurcation diagram of system \eref{eq03} for the equilibrium point $(0,C)$ with $(A,C,d)=(0.15,0.28,5)$ fixed. In the right panel, we show the real part of the dispersion relation $\lambda(k)$ as a function of the wave number squared for the system parameter $(A,C,d) = (0.15,0.28,5)$ fixed and $(\alpha)$ if $Q=0.5$ and $S=0.0701$ then $(0,C)$ is unstable with respect to small wave numbers and with respect to wave numbers near $\lambda_d^{\alpha}$; $(\beta)$ if $Q=0.55$ and $S=0.0701$ then $(0,C)$ is unstable only with respect to small wave numbers; $(\gamma)$ if $Q=0.575$ and $S=0.0701$ then the equilibrium point $(0,C)$ is stable in the ODE and PDE system.}
\label{F08}
\end{figure}

Next, we discuss the Turing conditions \eref{ODE} and \eref{PDE} for the equilibrium point $(0,C)$. The conditions in \eref{ODE} for the equilibrium point $(0,C)$ are met if the system parameters are such that $H_2/A-S<0$ and $-SH_2/A>0$. That latter implies that $H_2$ should be negative since the parameters $A$ and $S$ are positive. The first condition in \eref{PDE} for the equilibrium point $(0,C)$ is $dH_2>AS$, but $H_2<0$. Therefore, the conditions \eref{ODE} and \eref{PDE} for the equilibrium point $(0,C)$ cannot be met simultaneously, and thus we do not expect the formation of Turing patterns near $(0,C)$. The parameters space $(Q,S)$ of the equilibrium point $(0,C)$ is given in \Fref{F08} for the system parameters $(A,C)=(0.15,0.28)$ fixed. In the grey region ($\alpha$), the equilibrium point $(0,C)$ is unstable with respect to small wave numbers and with respect to wave numbers near $\lambda_d^{\alpha}$. In the green region ($\beta$), the equilibrium point $(0,C)$ is unstable only with respect to small wave numbers, while in the orange region ($\gamma$) the equilibrium point $(0,C)$ is stable in both the ODE and PDE system.

\subsubsection{Numerical Simulations of system \eref{eq03} near $(0,C)$}
Even though we do not expect Turing patterns, we present numerical solution to system \eref{eq03} in one and two-dimensional space for the system parameters $(A,C)=(0.15,0.28)$ fixed and initial condition near $(0,C)$. In particular, the initial condition is
\begin{equation}\label{icc}
u_0=0.012e^{-7x^2}\quad\text{and}\quad v_0=u_0+C=0.012e^{-7x^2}+0.28.
\end{equation} 
The numerical integration of system \eref{eq03} is performed under the same conditions used in \Sref{2Dp2}. 

\begin{figure}
\begin{center}
\includegraphics[width=7.5cm]{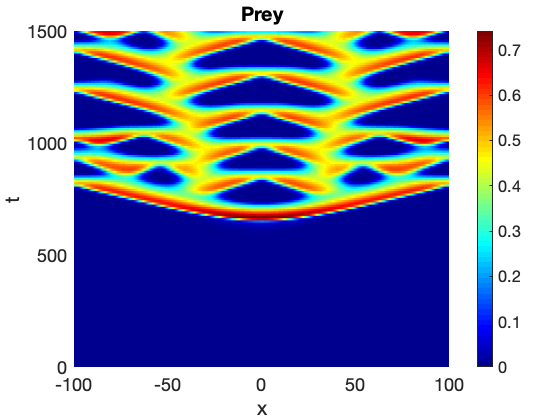}
\includegraphics[width=7.5cm]{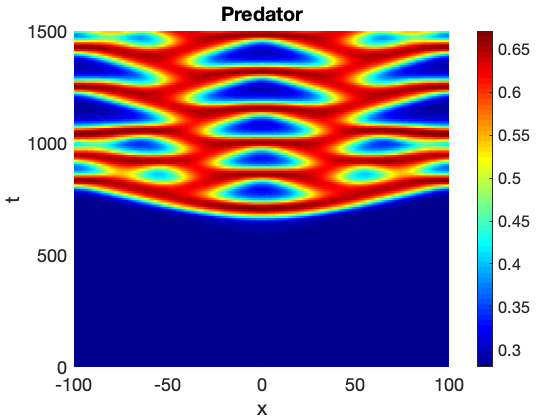}
\includegraphics[width=7.5cm]{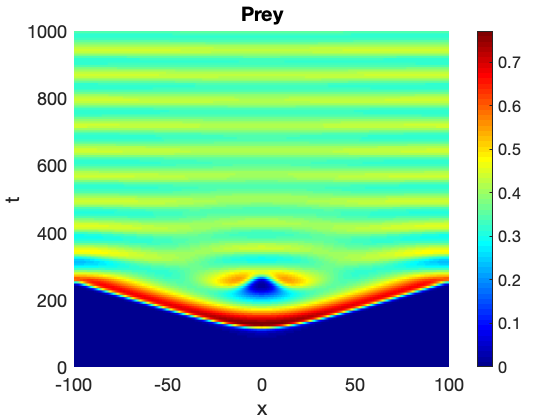}
\includegraphics[width=7.5cm]{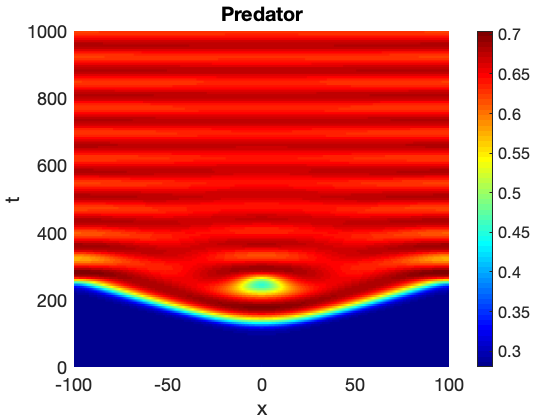}
\includegraphics[width=7.5cm]{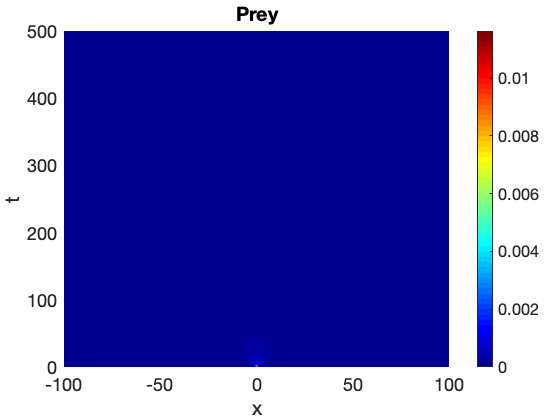}
\includegraphics[width=7.5cm]{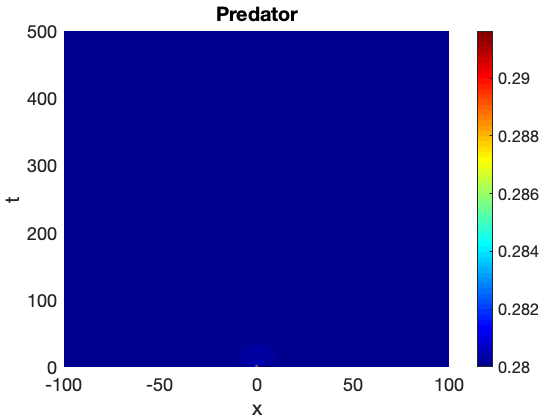}
\end{center}
\caption{Numerical simulation of system \eref{eq03} in one-dimensional space with system parameters $(A,C,d)=(0.15,0.28,5)$ fixed and initial condition defined in \eref{icc}. In the top panel $(Q,S)=(0.5,0.0701)$ and the equilibrium point $(0,C)=(0,0.28)$ is unstable with respect small wave numbers and wave numbers near $\lambda_d^{\alpha}$, see \Fref{F08}. We observe that the solution is irregular spatial pattern. In the middle panel $(Q,S)=(0.53,0.0701)$ and the equilibrium point $(0,C)=(0,0.28)$ is unstable with only respect small wave numbers. We observe that the solution is oscillatory in time with period $250(t)$. In the bottom panel $(Q,S)=(0.575,0.0701)$ and the equilibrium point $(0,C)=(0,0.28)$ is stable in the ODE and in the PDE. We observe that the initial condition evolves to $(0,C)$.} 
\label{F09}
\end{figure}

In region $(\alpha)$, see \Fref{F08}, the equilibrium point $(0,C)$ is unstable with respect to small wave numbers and with respect to wave numbers near $\lambda_d^{\alpha}$ and in $(\beta)$ the equilibrium point $(0,C)$ is unstable only with respect to small wave numbers. We observe that the initial condition \eref{icc} evolves to a spatial pattern that is oscillatory in time, see top panel of \Fref{F09}. In $(\gamma)$ the equilibrium point $(0,C)$ is stable in the ODE and PDE system. We observe that the initial condition \eref{icc} evolves, as expected, to the equilibrium point $(0,C)$, see bottom panel of \Fref{F09}. 

\begin{figure}
\begin{center}
\includegraphics[width=5.2cm]{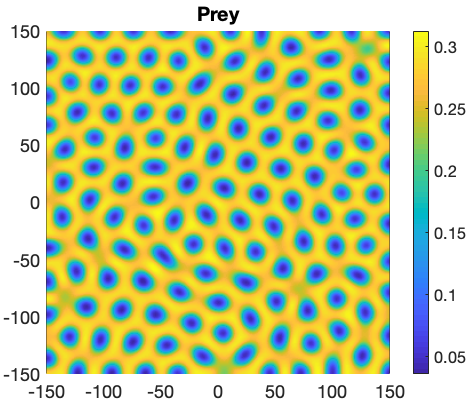}
\includegraphics[width=5.2cm]{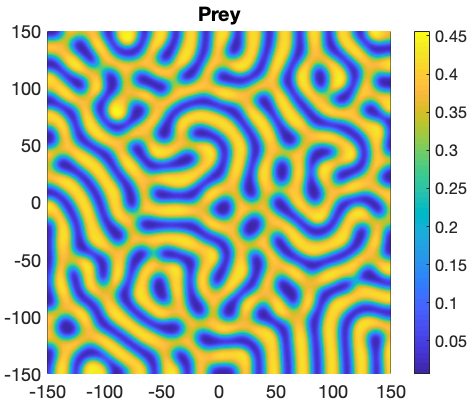}
\includegraphics[width=5.2cm]{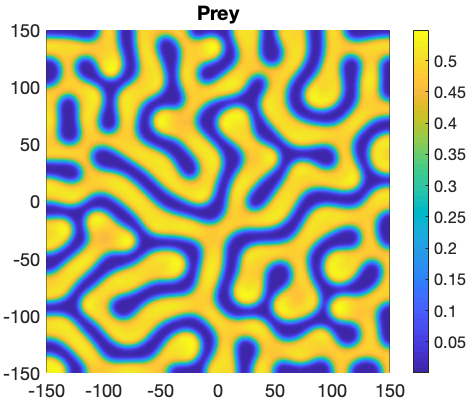}
\end{center}
\caption{Numerical simulation of the prey population for system \eref{eq03} in two-dimensional space for initial conditions near to the equilibrium point $(0,C)$. We consider system parameters $(A,C,Q,S)=(0.15,0.28,0.575,0.26)$ fixed and $d=5$ in the left panel; $d=9$ in the middle panel; and $d=55$ in the right panel. We observe the same type of spatial patterns as presented in \Fref{F07}. Note that these patterns are stationary in time.}
\label{F10}
\end{figure}

We also shortly discuss the spatio-temporal patterns of system \eref{eq03} in two-dimensional space for region $(\alpha)$ in \Fref{F08} where the equilibrium point $(0,C)$ is unstable with respect small wave numbers and wave numbers near $\lambda_d^{\alpha}$. The numerical integration of system \eref{eq03} is performed under the same conditions as in \Sref{2Dp2} and we consider the system parameter $(A,C,Q,S)=(0.15,0.28,0.575,0.26)$ fixed. The initial condition is a small random perturbation around the positive equilibrium point $(0,C)=(0,0.28)$. For the ratio of diffusivity $d=5$ we find cold-spot pattern over the whole domain, see left panel in  \Fref{F10}. When the ratio of diffusivity is being increased to $d=9$ we found that the cold-spot pattern started coalescing creating combination of labyrinthine and cold-spot pattern, see the middle panel in \Fref{F10}. Finally, by increasing the ratio of diffusivity up to $d=55$ we observe only labyrinthine pattern over the whole domain, see the right panel in \Fref{F10}. Note that all these patterns are stationary patterns as they remain unaltered with the further increase in time. Moreover, in \Fref{F10} we observe the same type of spatial pattern presented in \Fref{F07}, however the number of cold-spot pattern is approximately double than the number of cold-spot pattern presented in \Fref{F10} in the same domain. 

\section{Conclusion}\label{con}
In this manuscript, we study analytically and numerically the influence of diffusion on the pattern formation of the modified Holling--Tanner model \eref{eq01} with an alternative food source for the predator. We show that the pattern formation in the modified Holling-Tanner predator-prey model is rich and complex. In particular, we determine the temporal stability of the positive equilibrium point $P_2$ (see \Tref{T01}) and Turing instability of the same equilibrium point (see \Tref{T02}). We demonstrate the existence of a Turing instability of system \eref{eq03}, see region $(ii)$ in \Fref{F02} and \Fref{F03}. We show that in the neighbourhood of a Hopf bifurcation there exists a region in $(Q,S)$ parameter space (where $Q$ is related to the predation rate and $S$ is related to the intrinsic growth rate for the predator) in which the dynamics exhibits spatio-temporal behaviour that is influenced by a limit cycle. Moreover, we show that there exist conditions where the distribution of species oscillate in space and time, see region $(iii)$ in \Fref{F02} and top panel of \Fref{F05}. Furthermore, the numerical simulations in two-dimensional space of system \eref{eq03} show that, by increasing the diffusive constant $d$, the cold-spot pattern start to coalesce creating mixture Turing patterns (labyrinthine patterns and cold-spot patterns) and finally only labyrinthine patterns. The cold-spot patterns show that the prey population driven by the predator population leads to a lower proportion of the prey species in those regions. Then, by increasing the diffusive constant the population propagates in the space generating a population invasion, such that the holes connect each other forming tunnels with low population density, see \Fref{F07}. We also determine the temporal stability of the equilibrium point $(0,C)$. We found that the Turing conditions \eref{ODE} and \eref{PDE} for the equilibrium point $(0,C)$ cannot be met simultaneously and therefore we do not expect the formation of Turing patterns. We observe that initial conditions evolve to spatial patterns that are oscillatory in time when $(S,Q)$ are located in regions where a limit cycle exist, see \Fref{F08}.
\begin{figure}
\begin{center} 
\includegraphics[width=5.2cm]{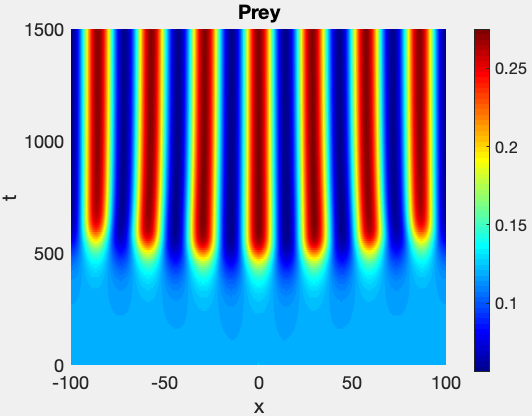}
\includegraphics[width=5.2cm]{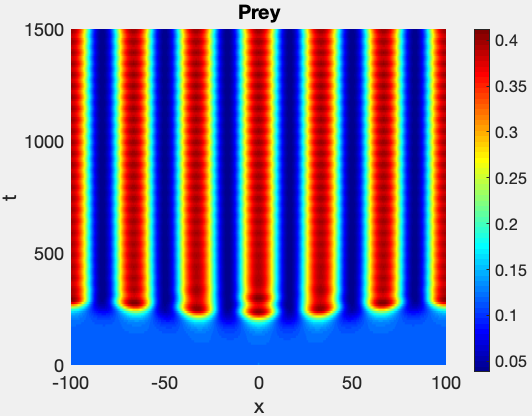}
\includegraphics[width=5.2cm]{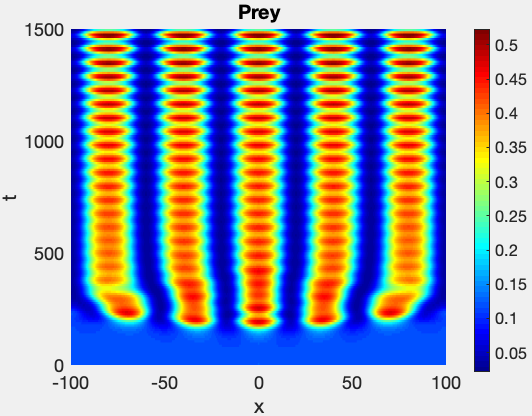}
\includegraphics[width=5.2cm]{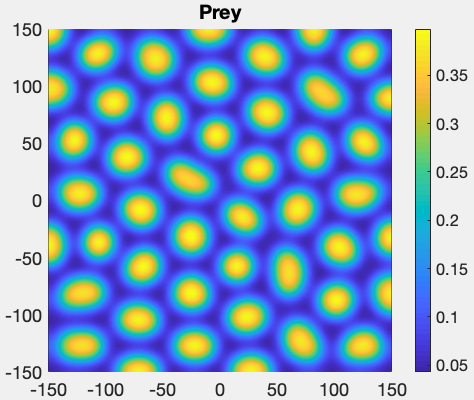}
\includegraphics[width=5.2cm]{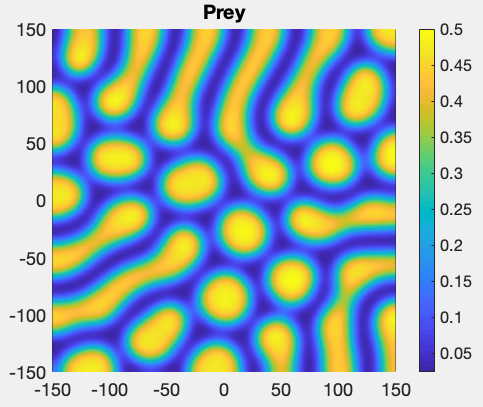}
\includegraphics[width=5.2cm]{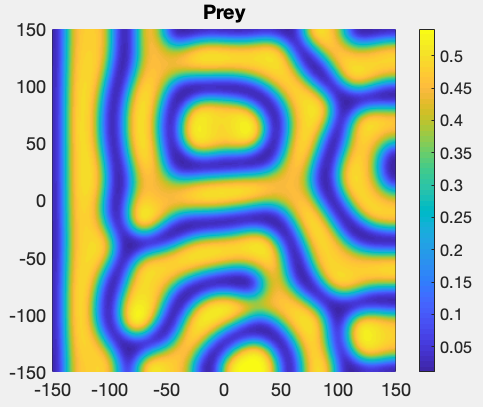}
\end{center}
\caption{In all cases we consider system parameters $(A,C,Q)=(0.15,0,2)$ fixed. On the top panel we show the spatial pattern of system \eref{eq03} in one-dimensional space with $t=1500$ and top left panel $S=0.28$, top middle panel $S=0.17$, and top right panel $S=0.13$. On the bottom panel we show the spatial pattern of system \eref{eq03} in two-dimensional space for the parameter $S=0.28$ also fixed and $d=15$ in the bottom left panel; $d=35$ in the bottom middle panel; and $d=65$ in the bottom right panel. Note that these Turing patterns are stationary in time.}
\label{F11}
\end{figure}

By taking the parameter $C=0$, that is, by removing the alternative food source for the predator system \eref{eq03} become singular for $u=0$. The temporal analysis of system \eref{eq03} with $C=0$ has been studied in \cite{saez}. The authors proved the existence of a non-hyperbolic equilibrium point $(0,0)$ and a saddle point $(1,0)$. In addition, there is always one equilibrium point in the first quadrant and this equilibrium point can be stable, unstable surrounded by a stable limit cycle, or stable surrounded by two limit cycles. The global stability of these periodic solutions was studied in \cite{lisena,hsu}. The spatio-temporal analysis of system \eref{eq03} with $C=0$ has been studied in \cite{chen,wang4}. The authors numerically showed that the model supports (Turing) patterns that are either periodic in space and stationary in time; or periodic in both in space and time, see top panel of \Fref{F11}. Moreover, the authors showed that system \eref{eq03} exhibits different Turing pattern formation in two-dimensional space such as hot-spot patterns, mixture Turing patterns (labyrinthine patterns and hot-spot patterns) and finally only labyrinthine patterns, see bottom panel of \Fref{F11}. We observe that in the model without alternative food the hot spot represent communities in which the prey, respectively the predator, interact (hot spot). While the modified model exhibit areas which are surrounded by communities (cold spot). 
\begin{figure}
\begin{center} 
\includegraphics[width=12cm]{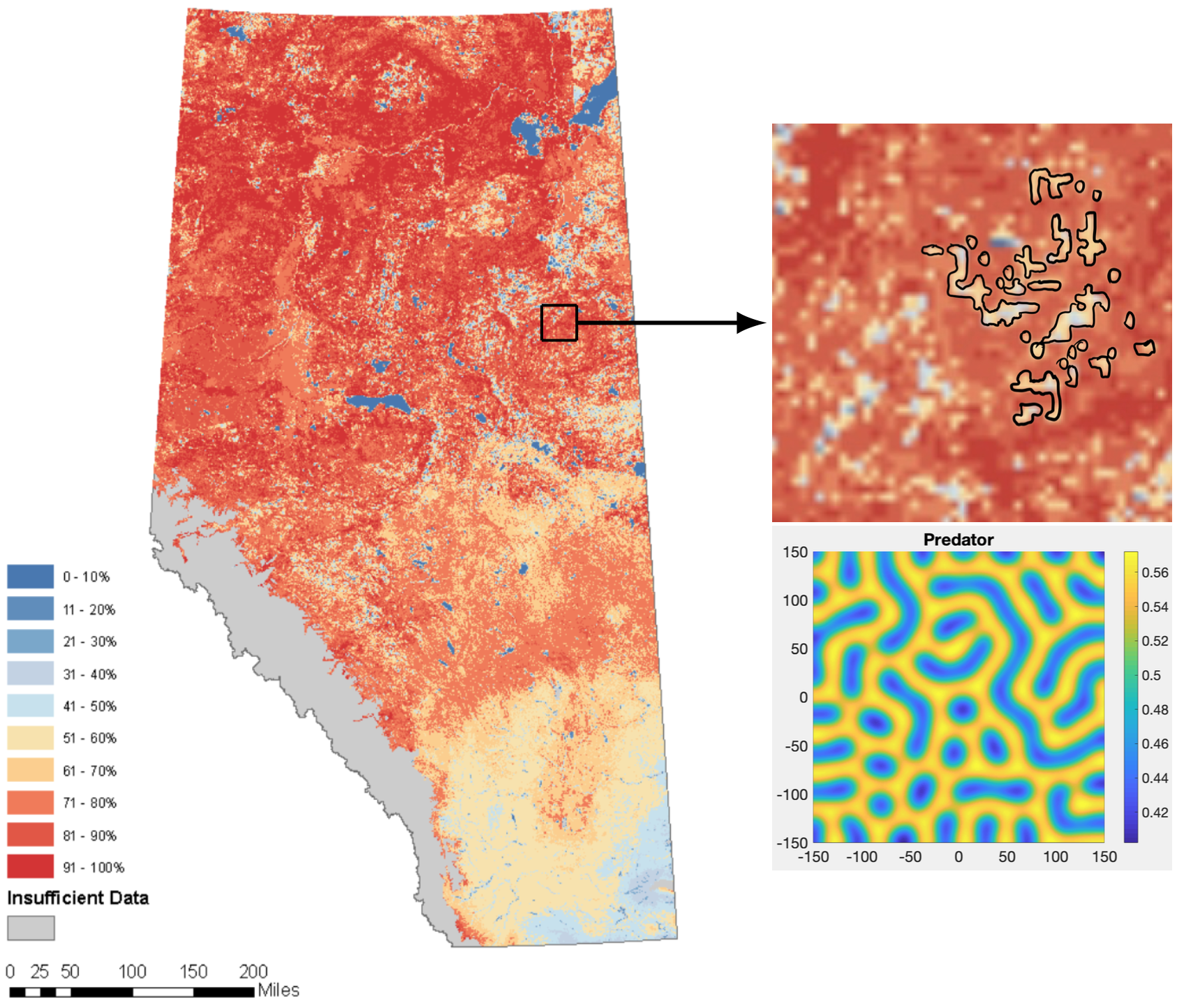}
\end{center}
\caption{The predicted relative abundance of the weasels and ermine in the forests of the Boreal Forest and Foothills Natural Regions. The map on the left is reproduced with permission from \cite{weasel}.}
\label{F12}
\end{figure}

The current relative abundance of the weasels observed in the Boreal Forest and Foothills Natural Regions shows that the distribution of this species increases, in general, from approximately 51\% up to 100\% of the abundance \cite{weasel}. In \Fref{F12} we observe that the distribution of weasels form the combination of two types of patterns, i.e. cold-spot patterns and labyrinthine patterns. This is similar to the Turing patterns presented in \Fref{F07}. Since the Turing patterns of the modified model showed in \Fref{F07} oscillate between $0.4$ and $0.55$, while the Turing pattern of the original model showed in \Fref{F11} oscillate between $0$ and $0.55$. Additionally, the type of patterns presented in the model without alternative food differ with the Turing pattern presented in the original diffusive Holling--Tanner model studied in \cite{chen,wang4}. In other words, the addition of the alternative food source for the predator in the model generate patterns which better represent the behaviour observed in the Boreal Forest and Foothills Natural Regions. Besides Turing patterns, numerical simulations also suggest that the spatial distribution of the weasels considering an alternative food source better represent the observed distribution of the population in this regions.  

\begin{figure}
\begin{center} 
\includegraphics[width=16cm]{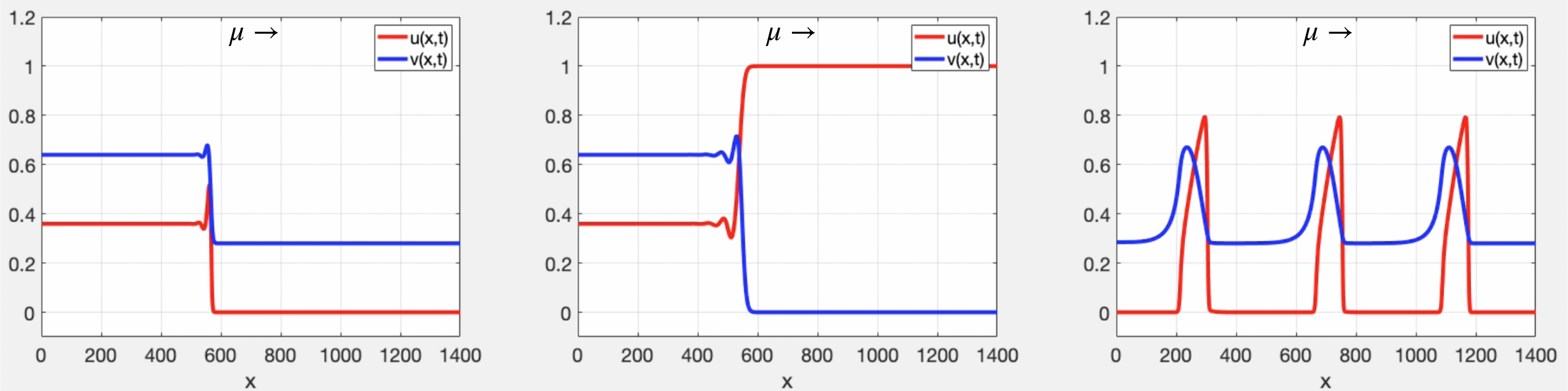}
\end{center}
\caption{Travelling wave solutions obtained by verifying the scaled intrinsic growth rate of the predator ($S$) and speed $\mu$. In the left and middle panel $S=0.18$ and in the right panel $S=0.03$. In the left panel, the travelling wave represents the invasion of the predator and the death of the prey population, while in the middle panel the travelling wave represents the extinction of predator and the stabilisation of the prey population.}
\label{F13}
\end{figure}

Numerical simulation indicate that the modified Holling--Tanner predator-prey model also supports travelling wave solutions. For the system parameter $(d,A,C,Q)=(5,0.15,0.28,0.51)$ fixed we numerically find that if $S=0.18$ (the scaled intrinsic growth rate of the predator) then there are travelling wave solutions connecting the equilibrium points $(0,C)=(0,0.28)$ and $(1,0)$ with $P_2=(0.36,0.64)$, see left and middle panel of \Fref{F13}. Additionally, by reducing the parameter $S=0.03$ we observe that the model also supports wave trains, see right panel of \Fref{F13}. These travelling wave solutions represents the invasion of the predator and the death of the prey population and also the extinction of predator and the stabilisation of the prey population. The analysis of these traveling waves solutions is left as further work.

\section*{Acknowledgments}
Weasels and Ermine data (2004-2013) from the Alberta Biodiversity Monitoring Institute was used \cite{weasel}, in part, to create \Fref{F12}. More information on the ABMI can be found at: \href{www.abmi.ca}{www.abmi.ca}.

\bibliographystyle{plain}
\bibliography{References.bib}

\end{document}